%% file: main.tex
\documentclass[11pt]{article}

\usepackage{amssymb}
\usepackage{amsmath}
\usepackage{setspace}
\usepackage{amsfonts}
\usepackage{algorithm,algorithmicx,algpseudocode}
\usepackage{array}
\usepackage[caption=false,font=normalsize,labelfont=sf,textfont=sf]{subfig}
\usepackage{textcomp}
\usepackage{stfloats}
\usepackage{url}
\usepackage{verbatim}
\usepackage{graphicx}
\usepackage{cite}
\usepackage{accents}
\usepackage{booktabs}
\usepackage{enumitem}
\usepackage{multirow}
\usepackage{siunitx}
\usepackage[dvipsnames]{xcolor}
\newlist{abbrv}{itemize}{1}
\setlist[abbrv,1]{label=,labelwidth=1.0in,align=parleft,itemsep=0.1\baselineskip,leftmargin=!}

\usepackage{pifont}

\usepackage{geometry}
\usepackage{tikz}
\usetikzlibrary{shapes.geometric}

\usepackage{balance}
\usepackage{todonotes}

\usepackage[dvipsnames]{xcolor}
\usepackage{amsmath}
\usepackage{multicol}
\usepackage{graphicx}
\usepackage{amsfonts,bm}
\usepackage{amsthm}
\usepackage{array}
\usepackage{textcomp}
\usepackage{stfloats}
\usepackage{url}
\usepackage{verbatim}
\usepackage{cite}
\usepackage{multirow}
\usepackage{booktabs}
\usepackage{hyperref}
\usepackage{stmaryrd}
\usepackage{amssymb}
\usepackage{pifont}
\usepackage{subcaption}
\usepackage{boldline}
\usepackage{cite}
\usepackage{enumitem}
\usepackage{calc}
\usepackage{multibib}
\usepackage{nomencl}
\usepackage{etoolbox}
\usepackage{mathtools}
\usepackage{cleveref}
\newcommand{\mmin}{\operatornamewithlimits{\text{\textbf{min}}}}

\DeclareFontFamily{U}{mathx}{}
\DeclareFontShape{U}{mathx}{m}{n}{<-> mathx10}{}
\DeclareSymbolFont{mathx}{U}{mathx}{m}{n}
\DeclareMathAccent{\widehat}{0}{mathx}{"70}
\DeclareMathAccent{\widecheck}{0}{mathx}{"71}

\newcommand{\pd}{\mathbf{d}}
\newcommand{\pfmax}{\mathbf{\bar{f}}}
\newcommand{\pfmin}{\mathbf{\ubar{f}}}
\newcommand{\pg}{\mathbf{p}}
\newcommand{\pf}{\mathbf{f}}
\newcommand{\pgmax}{\mathbf{\bar{p}}}
\newcommand{\res}{\mathbf{r}}
\newcommand{\resmax}{\mathbf{\bar{r}}}
\newcommand{\xith}{\mathbf{\xi}_{\text{th}}}

\newcommand{\xir}{\xi_{\text{r}}}

\newcommand{\cG}{\mathcal{G}} 
\newcommand{\cE}{\mathcal{E}} 
\newcommand{\cN}{\mathcal{N}} 
\newcommand{\cL}{\mathcal{L}} 
\newcommand{\slack}{\bm{\xi}} 
\newcommand{\slackapp}{\tilde{\bm{\xi}}}
\newcommand{\load}{\mathbf{d}} 
\newcommand{\g}{\mathbf{p}} 
\newcommand{\f}{\mathbf{f}} 
\newcommand{\gk}{\mathbf{p}_k} 
\newcommand{\rhok}{\bm{\rho}_k} 
\newcommand{\nk}{n_k} 

\newcommand{\Kg}{\mathcal{K}_g} 
\newcommand{\Ke}{\mathcal{K}_e}

\newcommand{\LODF}{\mathbf{L}} 
\newcommand{\flb}{\underline{\mathbf{f}}} 
\newcommand{\fub}{\overline{\mathbf{f}}} 
\newcommand{\glb}{\underline{\mathbf{p}}} 
\newcommand{\gub}{\overline{\mathbf{p}}} 
\newcommand{\slackpen}{M_\xi} 
\newcommand{\objc}{\mathbf{c}} 

\newcommand{\x}{\mathbf{x}}
\newcommand{\z}{\mathbf{z}}
\newcommand{\p}{\mathbf{p}}
\newcommand{\pnet}{P_{\theta}}
\newcommand{\dnet}{D_{\phi}}
\newcommand{\y}{\mathbf{y}}

\newcommand{\obj}{f_{\x}}
\newcommand{\obji}{f_{\x_i}}
\newcommand{\ineq}{\mathbf{g}_{\x}}
\newcommand{\eq}{\mathbf{h}_{\x}}
\newcommand{\eqi}{\mathbf{h}_{\x_i}}
\newcommand{\deq}{\bm{\lambda}}

\newcommand{\viol}[1]{\nu\left(#1\right)}
\newcommand{\bflam}{\mathbf{\lambda}}
\newcommand{\bfnu}{\mathbf{\nu}}

\newcommand{\rhomax}{\rho_{\text{max}}}

\newcommand{\norm}[1]{\left\lVert#1\right\rVert} 

\newcommand{\linfnorm}[1]{\left\lVert#1\right\rVert_\infty} 

\newcommand{\Mpb}{M_{\text{pb}}}
\newcommand{\Mres}{M_{\text{r}}}
\newcommand{\Mth}{M_{\text{th}}}

\newcommand{\ubar}[1]{\underline{#1}}

\DeclareMathOperator*{\argmin}{argmin}

\newcommand{\citep}[1]{\cite{#1}}
\usepackage{authblk}
\usepackage{float}
\floatstyle{ruled}
\newfloat{model}{thp}{lop}
\floatname{model}{Model}

\title{Optimization Learning}
\author{Pascal Van Hentenryck} 
\affil{NSF AI Institute for Advances in Optimization (AI4OPT) \\ Georgia Institute of Technology \\
Coda Building 12th/14th floor \\
756 W Peachtree St NW, Atlanta, GA 30308 \\
Email: pvh@gatech.edu}

\begin{document}
\maketitle

\begin{abstract}
This article introduces the concept of {\em optimization learning}, a
methodology to design optimization proxies that learn the input/output
mapping of parametric optimization problems. These optimization
proxies are trustworthy by design: they compute feasible solutions to
the underlying optimization problems, provide quality guarantees on
the returned solutions, and scale to large instances. Optimization
proxies are differentiable programs that combine traditional deep
learning technology with repair or completion layers to produce
feasible solutions. The article shows that optimization proxies can be
trained end-to-end in a self-supervised way. It presents methodologies
to provide performance guarantees and to scale optimization proxies to
large-scale optimization problems. The potential of optimization
proxies is highlighted through applications in power systems and, in
particular, real-time risk assessment and security-constrained optimal
power flow.
\end{abstract}

\input{notations}

\input{introduction}

\input{learning_task}

\input{formulation.tex}

\input{erm}

\input{primal_proxies}

\input{dual_proxies}

\input{primal_dual_proxies}

\input{conclusion}

\section*{Acknowledgements}

This research was partly supported by NSF award 2112533 and ARPA-E
PERFORM award AR0001136. Special thanks to the Optimization Proxy and
the Power Systems teams at AI4OPT and, in particular, Wenbo Chen,
Guancheng Qiu, Michael Klamkin, Terrence Mak, Seonho Park, and
Matthieu Tanneau for their invaluable contributions to this
research. Thanks to Ferdinando Fioretto and Andre Velloso for their
initial collaborations that shaped this research too.

\bibliographystyle{plain}
\bibliography{refs.bib}

\end{document}

%% file: notations.tex
\newcommand{\e}{\mathbf{e}}  
\newcommand{\bb}{\mathbf{b}}

\newcommand{\pgmin}{\mathbf{\ubar{p}}}
\newcommand{\rup}{R^{\uparrow}}
\newcommand{\rdn}{R^{\downarrow}}

\newcommand{\bA}{\mathbf{A}}
\newcommand{\bc}{\mathbf{c}}
\newcommand{\bx}{\mathbf{x}}
\newcommand{\by}{\mathbf{y}}
\newcommand{\bz}{\mathbf{z}}
\newcommand{\bw}{\mathbf{w}}
\newcommand{\bzl}{\mathbf{z}^{l}}
\newcommand{\bzu}{\mathbf{z}^{u}}
\newcommand{\bl}{\mathbf{l}}
\newcommand{\bu}{\mathbf{u}}
\newcommand{\be}{\mathbf{e}}
\newcommand{\bmu}{{\bm\mu}}
\newcommand{\bs}{\mathbf{s}}

%% file: introduction.tex
\section{Introduction}
\label{sec:intro}

Optimization technologies have been widely successful in industry:
they dispatch power grids, route transportation and logistic systems,
plan and operate supply chains, schedule manufacturing systems,
orchestrate evacuations and relief operations, clear markets for organ
exchanges to name only a few. Yet there remain applications where
optimization technologies still face computational challenges,
including when solutions must be produced in real time or when there
are planners and operators in the loop who expect fast interactions
with an underlying decision-support system.

Many engineering applications, however, operate on physical
infrastructures that change relatively slowly. As a result,
optimization technologies are often used to solve the same core
applications repeatedly on instances that are highly similar in
nature. In addition, it is reasonable to assume that these instances
followed a distribution learned from historical data, possibly
augmented with forecasts for future conditions. The overarching
scientific question is whether machine learning can help
optimization technologies expand the realm of tractable
optimization problems. 

This paper, which summarizes a significant part of my keynote at the
EURO-2024 conference in Copenhagen, reviews the concept of {\em
  optimization learning}, to address this challenge. The fundamental
idea behind optimization learning is to learn the input/output mapping
of optimization problems, i.e., the mapping from inputs to optimal
solutions. The paper presents three fundamental methodologies in
opimization learning: primal optimization proxies, dual optimization
proxies, and primal-dual learning. Primal optimization proxies return
(near-optimal) feasible solutions to parametric optimization problem
through a combination of deep learning and repair layers. Dual
optimization proxies return (near-optimal) dual solutions to
parametric optimization problems through a combination of deep
learning and completion layers. Primal-dual learning mimics
traditional augmented Lagrangian methods to provide primal proxies
that learn both a primal and a dual network.

These methodologies are illustrated on two significant applications in
power systems operations: the real-time risk assessment of a
transmission systems and the security-constrained optimal power under
N-1 Generator and Line Contingencies. In each case, optimization
learning brings orders of magnitude improvements in efficiency, making
it possible to solve the applications in real time with high accuracy,
an outcome that could not have been achieved by state-of-the-art
optimization technology.

The rest of this paper is organized as follows. Section
\ref{section:learning_task} presents the problem formulation, i.e.,
the learning task considered in the paper.  Section \ref{section:ED}
presents a running example for illustrating the methodologies: the
economic dispatch problem run every five minutes by independent system
operators in the US.  Section \ref{section:erm} discusses empirical
risk minimization in machine learning and baseline approaches for
approaching optimization learning. The section provides the background
for understanding the subsequent section and, in particular, a
high-level view of how machine learning models are trained.  Section
\ref{section:primal_proxies} presents primal optimization and their
application to real-time risk assessment. Section
\ref{section:dual_proxies} presents dual optimization proxies and their
application to DC optimal power flow. Section \ref{section:PDL}
presents primal-dual learning and its application to
security-constrained optimal power flows.  Section
\ref{sec:conclusion} concludes the paper and identifies future
research directions.

The paper unifies and reformulates, in a general framework, prior
research presented in several papers
\citep{E2ELR,klamkin2024dual,PDL,PDLSCOPF}, where extensive
literature reviews are available. ``Recent'' reviews of constrained
optimization learning and machine learning for optimal power flows are
available in \citep{ijcai2021p610,ML4OPFTutorial}. For space reasons,
this paper does not review the applications of optimization learning
for combinatorial optimization. Some applications in supply chains,
manufacturing, and transportation, where the optimization problems
contain discrete variables, can be found in
\citep{jss2022,ojha2024optimizationbased,DBLP:journals/jair/YuanCH22,DBLP:journals/corr/abs-2301-09703}.
For a general introduction to machine learning for combinatorial
optimization, see \citep{BENGIOLODI}.

%% file: learning_task.tex
\section{Problem Formulation: The Learning Task}
\label{section:learning_task}

This paper considers applications that require the solving of {\em
  parametric optimization problems}
\begin{equation}
\begin{aligned}
\label{eq:arg-opt}  
\Phi(\mathbf{x}) = \argmin_{\y} \obj(\y) \mbox{ subject to }  \eq(\y) = 0 \ \& \ \ineq(\y) \geq 0,
\end{aligned}
\end{equation}
where $\x$ represents instance parameters that determine the objective
function $\obj$ and the constraints $\eq$ and $\ineq$. Such parametric
optimization can be viewed as mapping from an input $\x$ to an output
$\y \Phi(\x)$ that represents its optimal solution (or a selected
optimal solution). In addition, for many applications, it is
reasonable to assume that the instance parameters $\x$ are
characterized by a distribution learned from historical data, possibly
augmented to cover future conditions with forecasting algorithms. The
goal of this paper is to study how to use the fusion of machine
learning and optimization for ``solving'' (unseen) optimization
instances $\Phi(\x_j)$ orders of magnitude faster than a
state-of-the-art optimization solver can. This task is called {\em
  optimization learning} in the rest of this paper.

More formally, given a probability distribution ${\cal P}$ of instance
parameters, the goal is to learn a {\em primal optimization proxy}
$\Phi^\uparrow: \Re^m \mapsto \Re^a$ that produces feasible solutions,
i.e.,
\[
\forall \x \sim {\cal P}: \eq(\Phi^\uparrow(\x)) = 0 \ \& \ \ineq(\Phi^\uparrow(\x)) \geq 0
\]
and a {\em dual optimization proxy} $\Phi^\downarrow: \Re^m \mapsto
\Re^b$ that returns valid lower bounds, i.e.,
\[
\forall \x \sim {\cal P} : \obj(\Phi^\downarrow(\x)) \leq \obj(\Phi(\x)).
\]

%% file: formulation.tex
\section{The Economic Dispatch Optimization}
\label{section:ED}

This paper illustrates optimization learning on a running example and
its variants: the Economic Dispatch (ED) optimization with reserve
requirements. The ED optimization is run every five minutes in the US
and is the backbone of power system operations.  The formulation
described in Model \ref{fig:ED} captures the essence of how real-time
markets are cleared by Independent Systems Operators (ISOs) in the
United States, although the reality is obviously more complex.

\begin{model}[!t]
    \caption{The Economic Dispatch Formulation.}
    \label{fig:ED}
    \begin{subequations}
    \label{eq:ED}
    \begin{flalign}
        \min_{\pg, \res, \xith} \quad & c(\pg) + \Mth \| \xith \|_{1} & & \\
        \text{s.t.} \quad
            & \mathbf{e}^{\top} \pg = \mathbf{e}^{\top} \pd, & & \label{eq:ED:power_balance}\\
            & \mathbf{e}^{\top} \res \geq R, & & \label{eq:ED:reserve_requirements}\\
            & \pg + \res \leq \pgmax, & &  \label{eq:ED:eco_max}\\
            & \mathbf{0} \leq \pg \leq \pgmax, & & \label{eq:ED:dispatch_bounds}\\
            & \mathbf{0} \leq \res \leq \resmax, & & \label{eq:ED:reserve_bounds}\\
            & \mathbf{\ubar{f}} -\xith \leq \mathbf{PTDF}(\pg - \pd)  \leq \mathbf{\bar{f}} +\xith, \label{eq:ED:PTDF} & & \\
            & \xith \geq \mathbf{0}. & & \label{eq:ED:thermal_slack_positive}
    \end{flalign}
    \end{subequations}
\end{model}

Constraints \eqref{eq:ED:power_balance} and
\eqref{eq:ED:reserve_requirements} are the power balance and minimum
reserve requirement constraints.  Constraints \eqref{eq:ED:eco_max}
ensure that the active power and reserves of each generator does not
violate their maximum capacities.  Constraints
\eqref{eq:ED:dispatch_bounds} and \eqref{eq:ED:reserve_bounds} enforce
the limits on each generator dispatch and reserves. Constraints
\eqref{eq:ED:PTDF} express the thermal limits on each limits using a
Power Transfer Distribution Factor (PTDF) representation, which is the
state of the art in industry
\citep{Ma2009_MISO_SCED,Holzer2022_MISO_SFT}. The thermal limits are
soft constraints, which is how they are modeled by US ISOs (e.g.,
\citep{Ma2009_MISO_SCED,BPM002_D}). Their violations incur a high
cost, modeled using slack variables $\xith$ which are penalized in the
objective.

%% file: erm.tex
\section{Empirical Risk Minimization}
\label{section:erm}

Machine learning seems to be an ideal approach to ``solve'' these
parametric optimization problems '' and replace optimization
all together. Indeed, an optimization problem can be viewed as mapping
from its input $\x$ to an optimal solution $\y$, and deep learning
networks are universal approximators \citep{HORNIK1989359}. This
section presents baseline machine learning approaches for optimization
learning.

\subsection{Supervised Learning}

A key feature of optimization learning is its ability to generate
correct labels for each input: it suffices to solve the optimization
problem \eqref{eq:arg-opt}. Hence, to train a machine learning model,
a data set ${\cal I} = \{ \x_i \}_{i \in [n]}$ of inputs sampled from
the distribution ${\cal P}$ can be augmented with their optimal
solutions to obtain a training data set ${\cal D} = \{
(\x_i,\y_i=\Phi(\x_i)) \}_{i \in [n]}$. Testing and validation data
sets can be obtained similarly.

\begin{figure}[!t]
\centering \includegraphics[width=.70\columnwidth]{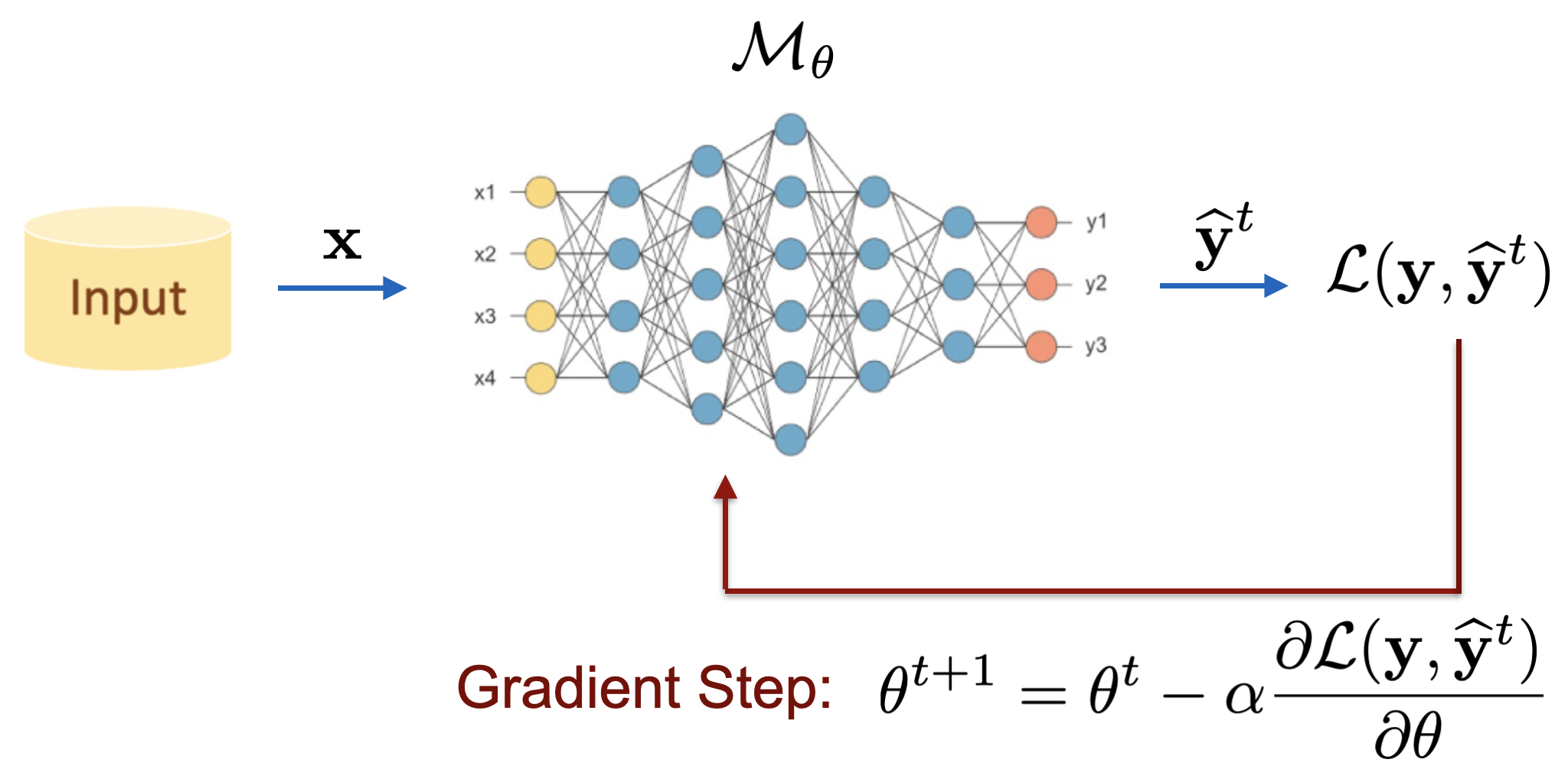}
\caption{The Training of Supervised Learning Models.}
\label{fig:training}
\end{figure}

{\em Supervised learning} can then used to fit the {\em learnable
  parameters} $\theta$ of a parametric machine learning model ${\cal M}_{\theta}$ by
minimizing a loss function of the form
\begin{equation}
\label{eq:erm}
\argmin_{\theta} \frac{1}{n} \sum_{i=1}^n {\cal L}(\y_i,{\cal M}_{\theta}(\x_i))
\end{equation}
where
\[
{\cal L}(\y,\widehat{\y}) = \| \y - \widehat{\y}) \|.
\]
The ``optimal'' value $\theta^*$ for the learnable parameters turns
the machine learning model into a function ${\cal M}_{\theta^*}$ that
approximates the parametric optimization $\Phi$. This minimization
problem is known as {\em an empirical risk minimization} and is
typically solved using stochastic gradient descent as illustrated in
Figure \ref{fig:training}. The machine learning training consists of
series of forward and backward passes. At iteration $t$, with
$\theta^t$ as the values of the learnable parameters, the forward pass
evaluates ${\cal M}_{\theta^t}(\x)$ to obtain $\widehat{\y}^t$, and the
backward pass updates the learnable parameters $\theta^t$ using a
gradient step
\[
\theta^{t+1} = \theta^t - \alpha \frac{\partial {\cal L}(\mathbf{y},\widehat{\mathbf{y}}^t)}{\partial \theta}.
\]

\noindent
Unfortunately, in engineering tasks, such an approximation is not a
primal optimization proxy: by virtue of being a regression, the machine learning
predictions are unlikely to satisfy the problem constraints, which may
have some significant consequences when optimization models are used
to plan, operate, or control physical infrastructures. In fact,
problem \eqref{eq:erm} should really be viewed as an {\em empirical
  minimization under constraints}, which raises significant challenges
for machine learning technologies.

\subsection{Lagrangian Duality in Supervised Learning}

Early research on this topic \citep{AAAI2020,MLKD2021} proposed to combine
machine learning and Lagrangian duality, using a parametric loss
function of the form
\[
{\cal L}^d_{\bflam,\bfnu}(\y,\widehat{\y}) = \| \y - \widehat{\y} \| + \bflam{} |\eq(\widehat{\y})| + \bfnu{} \max(0,-\ineq(\widehat{\y})).
\]
where $|\eq(\widehat{\y})|$ and $\max(0,-\ineq(\widehat{\y}))$ capture the
violations of the equality and inequality constraints. The multipliers
$\bflam{}$ and $\bfnu{}$ define the penalties for these constraint
violations; they can be trained by mimicking subgradient methods for
computing Lagrangian duals, alternating between training the machine
learning model and adjusting the multipliers with subgradient steps. At
each step $t$, the Lagrangian dual optimization solves a training
problem
\[
\theta^{t+1}  = \argmin_{\theta^t} \frac{1}{n} \sum_{i=1}^n {\cal L}^d_{\bflam^t,\bfnu^t}(\y_i,{\cal M}_{\theta^t}(\x_i))
\]
with the multipliers fixed at $\bflam^t$ and $\bfnu^t$. The multipliers are then
adjusted using the constraint violations, i.e.,
\begin{flalign*}
& \bflam^{t+1} = \bflam^t + \rho \frac{1}{n} \sum_{i=1}^n |\eq({\cal M}_{\theta^{t+1}}(\x_i))| \\
& \bfnu^{t+1}  = \bfnu^t + \rho \frac{1}{n} \sum_{i=1}^n \max(0,-\ineq({\cal M}_{\theta^{t<+1}}(\x_i))). 
\end{flalign*}
Experimentally, this Lagrangian dual approach has been shown to reduce
violations, sometimes substantially. However, being a regression, it
almost always violates constraints on unseen instances, and is not an
primal optimization proxy.

%% file: primal_proxies.tex
\section{Primal Optimization Proxies}
\label{section:primal_proxies}

\begin{figure}[!t]
\centering \includegraphics[width=.75\columnwidth]{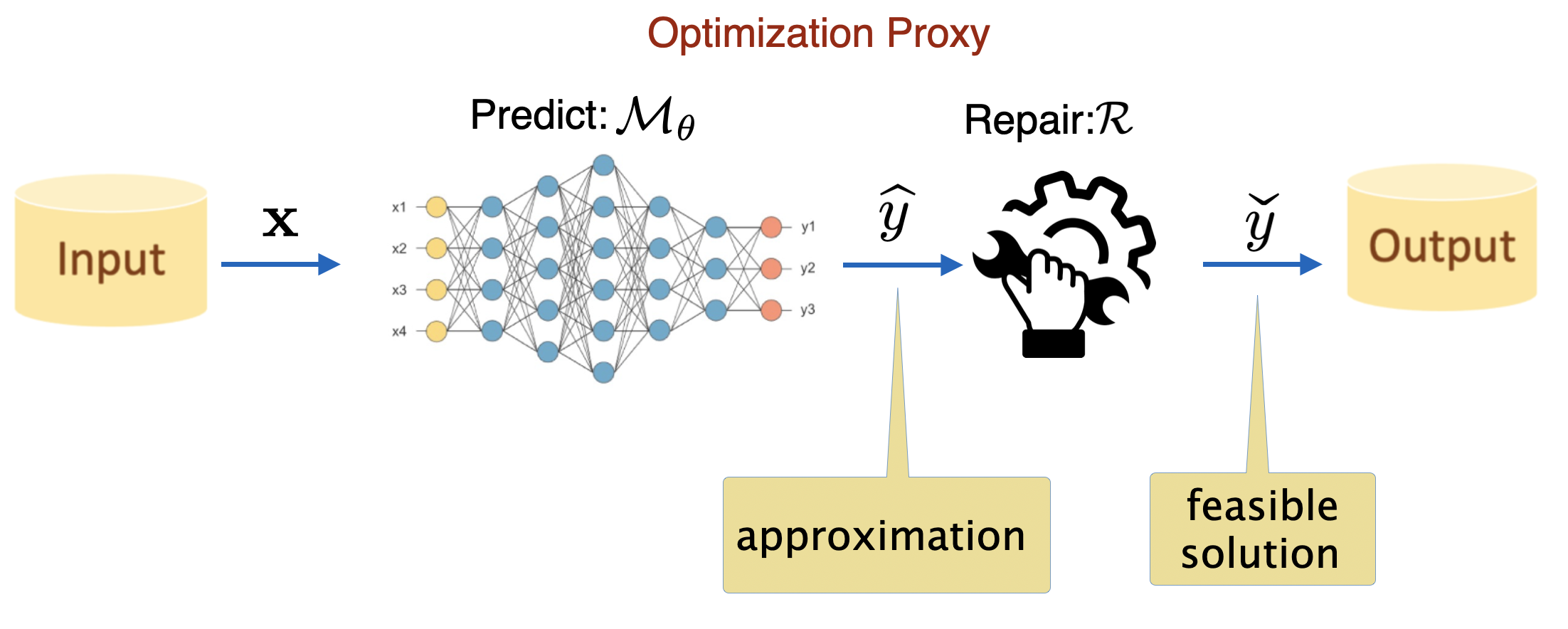}
\caption{The Architecture of Optimization Proxies.}
\label{fig:schematic}
\end{figure}

Optimization proxies were introduced to remedy the limitations of the
baseline machine learning models. Their high-level architecture is
depicted in Figure \ref{fig:schematic}. An optimization proxy
combines a parametric machine learning model ${\cal M}_{\theta}$
(typically a deep neural network) that produces an approximation
$\widehat{\y}$ with a {\em repair layer} ${\cal R}$ that transforms
$\widehat{\y}$ into a feasible solution $\widecheck{\y} = {\cal R}(\widehat{\y})$,
i.e., $\widecheck{\y}$ satisfies $\eq(\widecheck{\y}) = 0 \ \& \ \ineq(\widecheck{\y})
\geq 0$.  The repair layer can be thought of as a projection of
$\widehat{\y}$ into the feasible space of the optimization problem.

The simplest way to implement an optimization proxy is to train the
parametric machine learning as defined in Section \ref{section:erm} to
obtain $\theta^*$ and to use the repair layer at inference time to
obtain a feasible solution ${\cal R}({\cal M}_{\theta^*}(\x))$ for any
instance $\x$. The composition ${\cal R} \circ {\cal M}_{\theta^*}$ is a
primal optimization proxy. Moreover, an elegant way to implement the
repair layer is to solve the optimization model
\begin{equation}
{\cal R}(\widehat{\y}) = \argmin_{\y} \| \y - \widehat{\y} \|  \mbox{ subject to }  \eq(\y) = 0 \ \& \ \ineq(\y) \geq 0
\label{eq:projection}
\end{equation}
that finds the feasible point closest to $\widehat{\y}$. This is the
approach taken in \citep{velloso2021combining}.

\begin{figure}[!t]
\centering \includegraphics[width=.75\columnwidth]{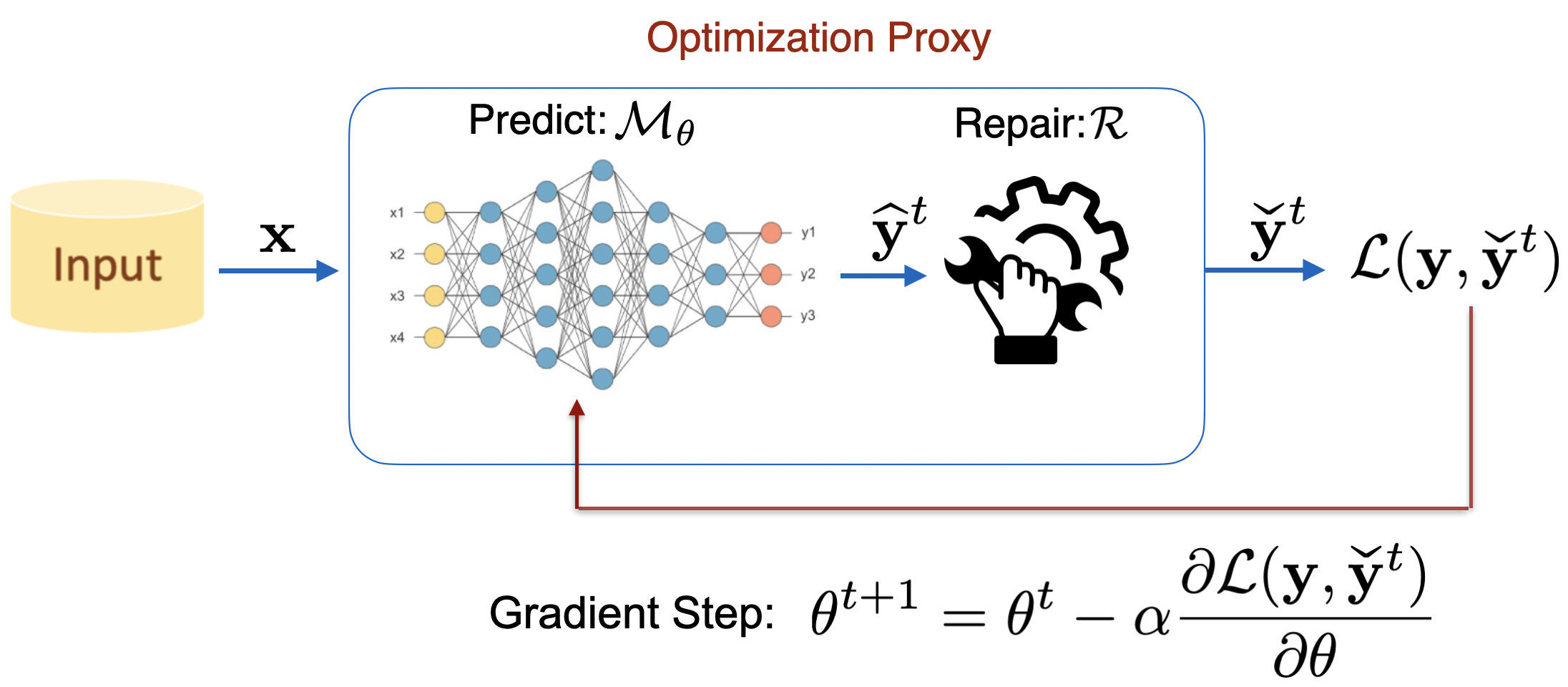}
\caption{End to End Training of Optimization Proxies.}
\label{fig:trainingproxy}
\end{figure}

An appealing alternative approach consists of training the
optimization proxy end-to-end as shown in Figure
\ref{fig:trainingproxy}. The challenge is to compute the gradient
\[
\frac{\partial {\cal L}(\mathbf{y},\widecheck{\mathbf{y}}^t)}{\partial \theta} = \frac{\partial {\cal L}(\mathbf{y},{\cal R}({\cal M}_{\theta^t}(\x)))}{\partial \theta}.
\]
By the chain rule, it can be computed as
\[
\frac{\partial {\cal L}(\mathbf{y},\widecheck{\mathbf{y}}^t)}{\partial \theta} =
\frac{\partial {\cal L}(\mathbf{y},\widecheck{\mathbf{y}}^t)}{\partial \widecheck{\mathbf{y}}^t} \;
\frac{\partial {\cal R}(\widehat{\mathbf{y}}^t)}{\partial \widehat{\mathbf{y}}^t} \;
\frac{\partial {\cal M}_{\theta^t}(\x)}{\partial \theta}   
\]
The difficulty, of course, is to differentiate the term $\frac{\partial
  {\cal R}(\widehat{\mathbf{y}}^t)}{\partial \widehat{\mathbf{y}}^t}$, which
is the topic of Sections \ref{section:implicit} and
\ref{section:dedicated}. 

\subsection{Implicit Differentiation}
\label{section:implicit}

A broad technique to differentiate through the repair layer is to use
Cauchy's implict function theorem \citep{ImplicitFunctionTheorem}.
Machine learning models can thus support two types of layers in their
networks: {\em explicit} layers that are specified by differentiable
functions and {\em implicit} layers whose outputs are the solutions to
systems of equations. Each layer should support the forward and
backward passes as presented in Section \ref{section:erm}.  Figure
\ref{fig:implicit} highlights the functionalities of explicit and
implicit layers. The forward pass of an implicit layer receives an
input $\x$ and finds an output $\y$ that satisfies a system of
equations $f(\x,\y) = 0$, i.e., it solves a system of equations during
training. The implicit function theorem is the key element for the
backward pass. If $g(\x)$ is a function that, given an input $\x$,
returns an output $\y$ that satisfies $f(\x,\y) = 0$. then the backward
pass should return $\frac{ \partial g(\x)}{\partial \x}$. In general,
there is no closed form for function $g$ but the implicit function
theorem gives a closed form for its gradient:
\[
\biggl[ \frac{\partial g_i(x)}{\partial x_j} \biggr] = \biggl[ J_{f,y}(x,g(x)) ] \biggr]^{-1} \biggl[ J_{f,x}(x,g(x)) ] \biggr]
\]
where
\[
J_{f,y}(a_1,a_2) = \biggl[ \frac{\partial f_i}{\partial y_i}(a_1,a2) \biggr]
\]
and similarly for $J_{f,x}$.

The implicit function theorem can be applied to the repair layer in
\eqref{eq:projection} by taking the KKT conditions of the projection
optimization. These define a system of equations; the forward pass of
this repair layer finds a solution to these KKT conditions; the
backward pass uses the implicit function theorem and applies the
closed form for the gradient of $g$. There are, of course, limitations
to this approach. On the one hand, solving the KKT conditions during
training for each instance and at inference time may be too costly. On
the other hand, for discrete problems, the KKT conditions are not
well-defined, but the primal proxy could use a continuous relaxation of
the optimization problem.

\begin{figure}[!t]
  \centering
  \begin{tabular}{|l|c|c|}
    \hline
    & Explicit Layer & Implicit Layer \\
    \hline
    Forward Pass &  $f(x)$  &  find $y$ such that  $f(x,y) = 0$ \\
    \hline
    Backward Pass & $\frac{\partial f(x)}{\partial x}$   &  $\frac{\partial g(x)}{\partial x}$  where $f(x,g(x))=0$\\
    \hline
  \end{tabular}
\caption{Explicit versus Implicit Layers in Deep Learning.}
\label{fig:implicit}
\end{figure}

Implicit layers have been a topic of interest in machine
learning. They are the foundations of what are called {\em deep
  declarative networks} and differentiable optimization layers (e.g.,
\citep{DeepDeclarativeNetworks,Amos_ICML2017,Agrawal2019}). Implicit
layers have also been heavily used in decision-focused learning, where
the goal is to train a machine learning model by taking into account
the effect of its predictions on a downstream optimization.  See, for
instance, \citep{Wilder_AAAI2019} and the survey in
\citep{ijcai2021p610} to see relationships between the various
approaches.

\subsection{Dedicated Repair Layers}
\label{section:dedicated}

\begin{figure}[!t]
\centering \includegraphics[width=.75\columnwidth]{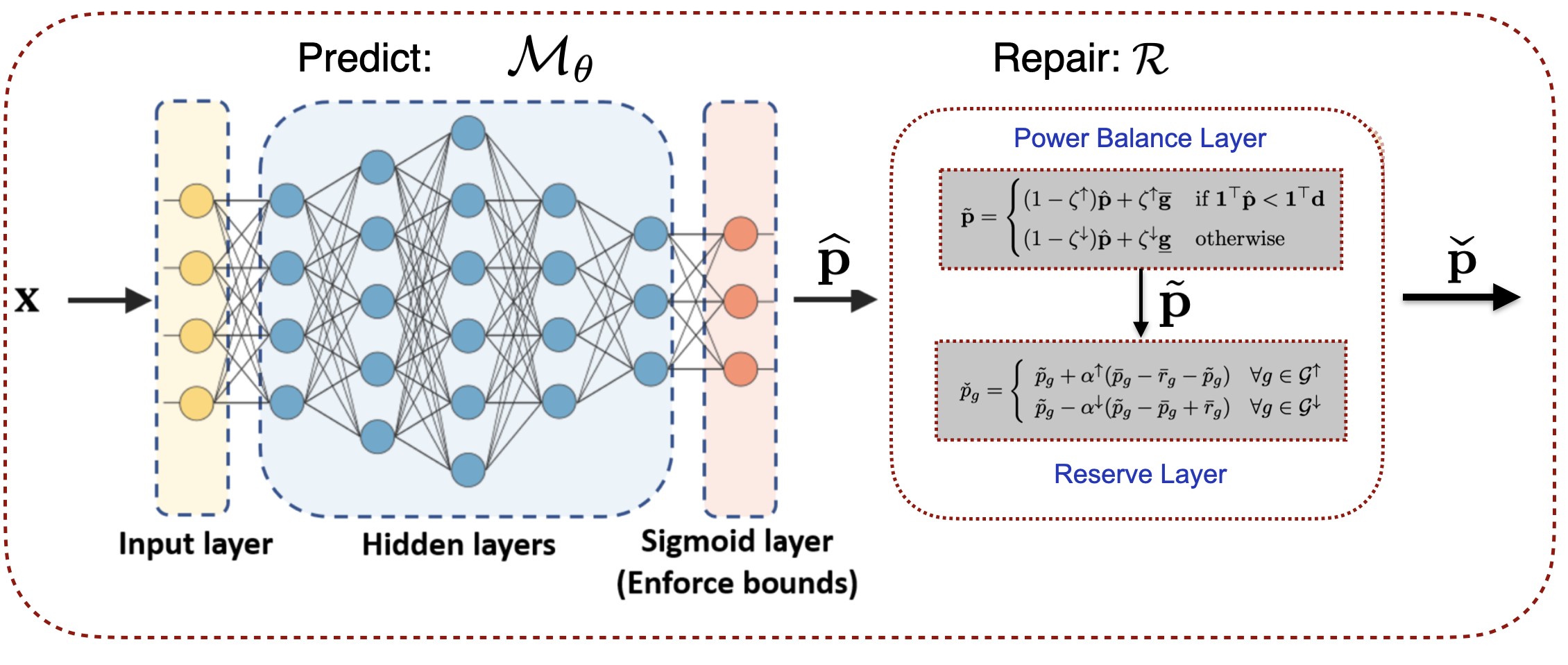}
\caption{The Optimization Proxy for the Economic Dispatch Problem (adapted from \citep{E2ELR}).}
\label{fig:EDE2E}
\end{figure}

Because of the high computational cost of implicit layers, it is
interesting to explore the design of dedicated repair layers and
learning architectures that ensure fast training and inference
times. This section illustrates this concept on the economic dispatch
optimization, where the challenge is to restore the feasibility of the
power balance and reserve constraints.  Figure \ref{fig:EDE2E} depicts
the architecture of a primal optimization proxy for the ED problem. It
uses a sigmoid layer to enforce the bound constraints, followed by a
repair layer for the power balance constraint and another repair layer
for the reserve constraints.

Consider first a dedicated repair layer for the balance constraint.
It takes as input an estimate $\widehat{p}_g$ for the power of each
generator $g$. It uses ideas from control
systems to scale all generators proportionally, up or down, to
transform $\widehat{\p} = {\cal M}_{\theta}(\x)$ into a solution $\tilde{\p}$
to the balance constraint as follows:
\begin{equation}
  \label{eq:frlayer}
\tilde{\p} = \begin{cases}
            (1-\zeta^\uparrow) \widehat{\p} + \zeta^\uparrow\gub     \;\;\;\;\text{if }\mathbf{1}^{\top}\widehat{\p}<\mathbf{1}^{\top}\load \\
            (1-\zeta^\downarrow) \widehat{\p} + \zeta^\downarrow\glb \;\;\;\;\text{otherwise}
\end{cases}
\end{equation}
where $\zeta^\uparrow$ and $\zeta^\downarrow$ are defined as,
\begin{align*}
    \zeta^\uparrow   = \frac{\mathbf{1}^{\top}\load-\mathbf{1}^{\top}\widehat{\p}}{\mathbf{1}^{\top}\gub-\mathbf{1}^{\top}\widehat{\p}} & &
    \zeta^\downarrow = \frac{\mathbf{1}^{\top}\widehat{\p}-\mathbf{1}^{\top}\load}{\mathbf{1}^{\top}\widehat{\p}-\mathbf{1}^{\top}\glb}
\end{align*}
and $\glb$ and $\gub$ are the lower and upper bounds on the
generators.  This repair layer \eqref{eq:frlayer} is differentiable
almost everywhere and can thus be naturally integrated in the training
of the machine learning model. Note that the formula for $\tilde{\p}$, i.e.,
\begin{align*}
\tilde{\p} = (1-\zeta^\uparrow) \widehat{\p} + \zeta^\uparrow\gub & & \mbox{ or } & & \tilde{\p} = (1-\zeta^\downarrow) \widehat{\p} + \zeta^\downarrow\glb
\end{align*}
can be generated during the forward pass based on the result of the
test $\mathbf{1}^{\top}\widehat{\p}<\mathbf{1}^{\top}\load$. It also
admits subgradients everywhere and can thus be backpropagated during
the backward pass. This also highlights the power of Dynamic
Computational Graphs (DCGs) and, more generally, {\em differentiable
  programming}: the computational graph does not have to be defined
statically and is not necessarily the same for each learning
iteration.

The reserve layer follows a similar approach but it is more involved.
During the forward phase, the proxy first partitions the generators into
(1) the set $\mathcal{G}^{\downarrow}$ of units whose power output can be decreased to
provide more reserves; and (2) the set $\mathcal{G}^{\uparrow}$ of units whose power output
can be increased with affecting their reserves
  \begin{align*}
    \mathcal{G}^{\uparrow} \gets \left\{ g \ \middle| \ \tilde{p}_{g} \leq \bar{p}_{g} - \bar{r}_{g} \right\}  & & 
    \mathcal{G}^{\downarrow} \gets \left\{ g \ \middle| \  \tilde{p}_{g} > \bar{p}_{g} - \bar{r}_{g} \right\}.
  \end{align*}
The amount $\Delta^{\downarrow}$ of additional reserves and additional power $\Delta^{\uparrow}$
that can be provided are given by
  \begin{align*}
    \Delta^{\downarrow} \gets \sum_{g \in \mathcal{G}^{\downarrow}} \tilde{p}_{g} - (\bar{p}_{g} - \bar{r}_{g}) & &
    \Delta^{\uparrow} \gets \sum_{g \in \mathcal{G}^{\uparrow}} (\bar{p}_{g} - \bar{r}_{g}) - \tilde{p}_{g} 
  \end{align*}
It is then possible to compute the additional reserve amount $\Delta$ that is available witout violating the
balance constraint
\[
\Delta \gets \max(0, \min(\Delta_{R}, \Delta^{\uparrow}, \Delta^{\downarrow}))
\]
where 
\[
\Delta_{R} \gets R - \sum_{g} \min \{ \bar{r}_{g}, \bar{p}_{g} - \tilde{p}_{g} \}
\]
denotes the shortage in reserves. The generators can then scaled proportionally by 
\begin{align}
\label{eq:reserves}
\widecheck{p}_{g} &= \left\{
                      \begin{array}{ll}
                            \tilde{p}_{g} + \alpha^{\uparrow} (\bar{p}_{g} - \bar{r}_{g} - \tilde{p}_{g}) & \forall g \in \mathcal{G}^{\uparrow} \\
                            \tilde{p}_{g} - \alpha^{\downarrow} (\tilde{p}_{g} - \bar{p}_{g} + \bar{r}_{g}) & \forall g \in \mathcal{G}^{\downarrow}
                        \end{array}
                        \right.
\end{align}
where $\alpha^{\uparrow} \gets \Delta / \Delta^{\uparrow}$ and
$\alpha^{\downarrow} \gets \Delta / \Delta^{\downarrow}$. Again, the
definition of $\widecheck{p}_{g}$ for each generator is differentiable.
Moreover, these two layers are guaranteed to find a feasible solution if one
exists. 

The overall architecture (see Figure \ref{fig:EDE2E}), called E2ELR,
is a differentiable program: it can be trained end to end efficiently
and runs in milliseconds during inference due to these dedicated
repair layers.

\subsection{Self-Supervised Learning}

Primal optimization proxies also offer an intriguing alternative to
supervised learning. Indeed, the approaches presented so far rely on
the availability of the data set ${\cal D} = \{ (\x_i,\y_i=\Phi(\x_i))
\}_{i \in [n]}$. However, since the parametric optimization problem is
available in explicit form, it is possible to define a {\em
  self-supervised} version of the learning task that only uses the
data set ${\cal I}$. Given a parametric primal optimization proxy
${\Phi}^{\uparrow}_{\theta}$, the self-supervised learning task
amounts to solving the optimization problem
\[
\min_{\theta} \frac{1}{n} \sum_{i=1}^n \obji(\Phi^{\uparrow}_{\theta}(\x_i))
\]
which is expressed in terms of the objective function $\obj$ of the
original parametric optimization problem. Note also that the learning
talk is different from the original optimization problem
\eqref{eq:arg-opt}: the optimization is over the learnable parameters,
not the decision variables of \eqref{eq:arg-opt}.  

For instance, E2ELR can be trained end-to-end using
self-supervised learning and the loss function
\[
c(\widecheck{\pg}) + \Mth \| \widecheck{\xi}_{\text{th}} \|_{1}
\]
which is then backpropagated through the repair layers to adjust the
parameters $\theta$ of the neural network. The training of a
self-supervised E2ELR does not require any solved instance.

This self-supervised approach has a significant benefit: it does not
rely on the availability of the optimal solutions for a large number
of instances, which may be costly to obtain in practice. In addition,
the objective function in supervised learning approaches may not be
perfectly aligned with the objective of the original optimization
problem, reducing the accuracy of the learning step. It should be
noted that supervised learning can also be used, say with the LD
approach proposed earlier. However, this provides less guidance to the
learning system than the labeled data in general, and the machine
learning model typically trade increased violations for better
objectives.  Early papers on self-supervised learning for optimization
can be found in (e.g., \citep{DC3,PDL,E2ELR,PDL,Huang2021_DeepOPF-NGT,wang2022fast}).

\subsection{Experimental Results}
\label{section:experiments}

\newcommand{\ieeeSmall}{\texttt{ieee300}}
\newcommand{\pegaseSmall}{\texttt{pegase1k}}
\newcommand{\rteLarge}{\texttt{rte6470}}
\newcommand{\pegaseLarge}{\texttt{pegase9k}}
\newcommand{\pegaseXLarge}{\texttt{pegase13k}}
\newcommand{\gocXXL}{\texttt{goc30k}}

This section demonstrates the capabilities of primal optimization
proxies, and E2ELR in particular, on power grids with up to 30,000
buses. All details can be found in \citep{Chen2023_E2ELR_arxiv}.

\subsubsection{Data Generation}
\label{sec:experiment:data}

    \begin{table}[!t]
        \centering
        \begin{tabular}{lcrrrrr}
            \toprule
            System & Ref
                & \multicolumn{1}{c}{$|\mathcal{N}|$}
                & \multicolumn{1}{c}{$|\mathcal{E}|$}
                & \multicolumn{1}{c}{$|\mathcal{G}|$} 
                & \multicolumn{1}{c}{${D_{\text{ref}}}^{\dagger}$}
                & \multicolumn{1}{c}{$\alpha_{\text{r}}$}
                \\
            \midrule
            \ieeeSmall      & \citep{UW_PowerSystemArchive} & 300   &	411 &	69  & 23.53    & 34.16\% \\
            \pegaseSmall    & \citep{Josz2016_ACOPF_PegaseRTE} & 1354  & 1991  & 260   & 73.06    & 19.82\% \\
            \rteLarge       & \citep{Josz2016_ACOPF_PegaseRTE} & 6470  & 9005  & 761   & 96.59    & 14.25\% \\
            \pegaseLarge    & \citep{Josz2016_ACOPF_PegaseRTE} & 9241  & 16049 & 1445  & 312.35   &  4.70\% \\
            \pegaseXLarge   & \citep{Josz2016_ACOPF_PegaseRTE} & 13659 & 20467 & 4092  & 381.43   &  1.32\% \\
            \gocXXL         & \citep{GoCompetition} & 30000 & 35393 & 3526  & 117.74   &  4.68\% \\
            \bottomrule
        \end{tabular}\\
        \caption{Selected test cases from PGLib \citep{PGLib}}
        \label{tab:PGLib}
    \end{table}

To obtain reasonable training, validation, and testing data sets,
benchmarks from the PGLib \citep{PGLib} library (v21.07) were modified
by perturbing the loads. Denote by $\pd^{\text{ref}}$ the nodal load
vector from the benchmarks. Load $i$ becomes
$\pd^{(i)} \, {=} \, \gamma^{(i)} \times \eta^{(i)} \times
\pd^{\text{ref}}$, where $\gamma^{(i)} \, {\in} \, \mathbb{R}$ is a
global scaling factor, $\eta \, {\in} \, \mathbb{R}^{|\mathcal{N}|}$
denotes load-level multiplicative white noise, and the multiplications
are element-wise.  $\gamma$ is sampled from a uniform distribution
$U[0.8, 1.2]$, and, for each load, $\eta$ is sampled from a log-normal
distribution with mean $1$ and standard deviation $5\%$. The PGLib
library does not include reserve information: for this reason, the
instance generation mimics the contingency reserve requirements used
in industry. It assumes $\bar{r}_{g} = \alpha_{\text{r}} \bar{p}_{g}$,
where $\alpha_{\text{r}} = 5 \times \| \bar{\pg} \|_{\infty} \times \|
\bar{\pg} \|_{1}^{-1}, $ ensuring a total reserve capacity 5 times
larger than the largest generator.  The reserve requirements of each
instance is sampled uniformly between 100\% and 200\% of the size of
the largest generator. Table \ref{tab:PGLib} presents the resulting
systems: it reports the number of buses ($|\mathcal{N}|$), branches
($|\mathcal{E}|$), and generators ($|\mathcal{G}|$), the total active
power demand in the reference PGLib case ($D_{\text{ref}}$, in GW),
and the value of $\alpha_{\text{r}}$. Large systems have a smaller
value of $\alpha_{\text{r}}$ because they contain significantly more
generators, but the size of the largest generator typically remains
roughly the same.  For every test case, $50,000$ instances are
generated and solved using Gurobi.  This dataset is then split into
training, validation, and test sets which comprise $40000$, $5000$,
and $5000$ instances.

\subsubsection{Baseline Models and Performance Metrics}
\label{sec:experiment:baselines}

The experiments compare E2ELR with four baseline models. The simplest
is a na\"ive, fully-connected DNN model without any feasibility layer
(DNN): it simply includes a Sigmoid activation layer to enforce the
generation bounds. The second baseline is a fully connected DNN model
with the DeepOPF architecture \citep{pan2020deepopf} (DeepOPF).  It
uses an equality completion to ensure the satisfaction of equality
constraints, but its output may violate inequality constraints.  The
third baseline is a fully connected DNN model with the DC3
architecture \citep{DC3} (DC3). This architecture uses a
partial repair layer, i.e., a fixed-step unrolled gradient descent to
minimize constraint violations. This layer mimics the implicit layer
described previously, but only performs a fixed number of gradient
steps for efficiency reasons. DC3 is not guaranteed to reach zero
violations however. The last model is a fully connected DNN model,
combined with the LOOP-LC architecture from
\citep{Li2022_LOOP-LC}. LOOP uses gauge functions to define a
one-to-one mapping between the unit hypercube and the set of feasible
solutions. The approach applies when all constraints are convex, the
feasible set is bounded, and a strictly feasible point is available
for each instance. LOOP was only evaluated on ED instances with no
reserves, since its gauge mapping does not support the reserve
constraint formulation in terms of generation variables only. The
baselines are not intended to be exhaustive, but rather to demonstrate
the potential of primal optimization proxies. Other approaches are
possible (e.g., \citep{Li2022_LOOP-LC,Tordesillas2023_rayen,konstantinov2024_NNHardConstraints}).

\subsubsection{Optimality Gaps}
\label{sec:results:gaps}

The first results measure the quality of the machine learning models using
the optimality gap which is generalized to include penalties for
constraint violations for a fair comparison.  Given an instance $\x$
with optimal solution $\pg^{*}$ and a prediction $\widehat{\pg}$, the
optimality gap is defined as $ \text{gap} = (\widehat{Z} - Z^{*}) \times |
Z^{*} |^{-1}, $ where $Z^{*}$ is the optimal objective value, and
$\widehat{Z}$ is the objective value of the prediction defined as
\begin{align}
        \label{eq:experiment:objective_penalized}
        c(\widehat{\pg}) + \Mth \| \widehat{\xith} \|_{1} + \Mpb |\mathbf{e}^{\top}(\widehat{\pg} - \pd)| + \Mres \xir(\widehat{\pg}).
\end{align}
The paper uses realistic penalty prices, based on the values used in
MISO operations
\citep{MISO_Schedule28_ReserveDemandCurve,MISO_Schedule28A_TranmissionDemandCurve}.

Table \ref{tab:res:opt_gap} reports the optimality gaps for the ED
instances and for the same instances with no reserve constraints to
evaluate LOOP (ED-NR). Both Supervised Learning (SL) and
Self-Supervised Learning (SSL) are considered.  Bold entries denote
the best performing method. The results show {\em E2ELR has the
  smallest optimality gaps across all settings}, primarily because it
always returns feasible solutions and avoids large penalties. The
performance of LOOP on the ED-NR is not strong, although it produces
feasible solutions: its non-convex gauge mapping guides the learning
to a local optimum after a few epochs. Statistics on power balance
violations for DNN, DeepOPF, and DC3 are reported in Table
\ref{tab:res:violation}. E2ELR with 
self-supervised training is particularly effective.  With the exception of \rteLarge{} (the real
French transmission grid), its performance improves with system sizes.

\begin{table}[!t]
  \small
        \centering
        \input{gaps_all.tex}

                \caption{Mean Optimality Gaps (\%) (from \citep{E2ELR}).}
        \label{tab:res:opt_gap}
    \end{table}

    \begin{table}[!t]
      \centering
      \small
            \input{viol_powerBalance_all}
        \caption{Violations of the Power Balance Constraint (from \citep{E2ELR}).}
        \label{tab:res:violation}
    \end{table}
        
\subsubsection{Computing Times}
\label{sec:results:times}

The second set of results evaluates the computational efficiency of
each ML model.  The computational efficiency is measured by (i) the
training time of ML models, including the data generation time when
applicable, and (ii) the inference time.  Note that ML models evaluate
\emph{batches} of instances, and inference times are reported per
batch of 256 instances.  Unless specified otherwise, average computing
times are arithmetic means; other averages use shifted geometric means
with a shift $s$ of 1\% for optimality gaps, and 1 p.u. for
constraint violations.

Table \ref{tab:exp:ED-R:train_time} reports the sampling and training
times.  The sampling times represent the total time to obtain the
optimal solutions to all instances (using the Gurobi solver on a
single thread) The training times for SL and SSL are comparable but
SSL has the significant benefit of not requiring optimal solutions to
be produced offline. The training time of DC3 is significantly higher
than other baselines because of its unrolled gradient steps: this is
the issue mentioned earlier on the cost of using complex implicit
layers.  The main takeaway is that {\em the self-supervised E2ELR
  needs less than an hour of training to achieve optimality gaps under
  $0.5\%$ for grids with thousands of buses.}

    \begin{table}[!t]
      \centering
        \small
        \begin{tabular}{llrrrrr}
            \toprule
            Loss & System & \multicolumn{1}{c}{Sampling} & \multicolumn{1}{c}{DNN} & \multicolumn{1}{c}{E2ELR} & \multicolumn{1}{c}{DeepOPF} & \multicolumn{1}{c}{DC3} \\
            \midrule
            SL
            & \ieeeSmall       & \qty[mode=text]{ 0.2}{\hour} & \qty[mode=text]{12}{\minute}  & \qty[mode=text]{43}{\minute} & \qty[mode=text]{43}{\minute} & \qty[mode=text]{115}{\minute} \\
            & \pegaseSmall     & \qty[mode=text]{ 0.8}{\hour} & \qty[mode=text]{14}{\minute}  & \qty[mode=text]{19}{\minute} & \qty[mode=text]{19}{\minute} & \qty[mode=text]{53 }{\minute} \\
            & \rteLarge        & \qty[mode=text]{ 4.6}{\hour} & \qty[mode=text]{14}{\minute}  & \qty[mode=text]{19}{\minute} & \qty[mode=text]{19}{\minute} & \qty[mode=text]{71 }{\minute} \\
            & \pegaseLarge     & \qty[mode=text]{14.0}{\hour} & \qty[mode=text]{15}{\minute}  & \qty[mode=text]{22}{\minute} & \qty[mode=text]{22}{\minute} & \qty[mode=text]{123}{\minute} \\
            & \pegaseXLarge    & \qty[mode=text]{22.7}{\hour} & \qty[mode=text]{16}{\minute}  & \qty[mode=text]{27}{\minute} & \qty[mode=text]{27}{\minute} & \qty[mode=text]{126}{\minute} \\
            & \gocXXL          & \qty[mode=text]{65.9}{\hour} & \qty[mode=text]{32}{\minute}  & \qty[mode=text]{39}{\minute} & \qty[mode=text]{38}{\minute} & \qty[mode=text]{129}{\minute} \\
            \midrule
            SSL
            & \ieeeSmall       & -- & \qty[mode=text]{21}{\minute}  & \qty[mode=text]{37}{\minute} & \qty[mode=text]{37}{\minute} & \qty[mode=text]{131}{\minute} \\
            & \pegaseSmall     & -- & \qty[mode=text]{ 6}{\minute}  & \qty[mode=text]{19}{\minute} & \qty[mode=text]{19}{\minute} & \qty[mode=text]{ 67}{\minute} \\
            & \rteLarge        & -- & \qty[mode=text]{12}{\minute}  & \qty[mode=text]{21}{\minute} & \qty[mode=text]{21}{\minute} & \qty[mode=text]{ 71}{\minute} \\
            & \pegaseLarge     & -- & \qty[mode=text]{20}{\minute}  & \qty[mode=text]{24}{\minute} & \qty[mode=text]{24}{\minute} & \qty[mode=text]{123}{\minute} \\
            & \pegaseXLarge    & -- & \qty[mode=text]{13}{\minute}  & \qty[mode=text]{22}{\minute} & \qty[mode=text]{22}{\minute} & \qty[mode=text]{125}{\minute} \\
            & \gocXXL          & -- & \qty[mode=text]{52}{\minute}  & \qty[mode=text]{53}{\minute} & \qty[mode=text]{53}{\minute} & \qty[mode=text]{128}{\minute} \\
            \bottomrule
        \end{tabular}
        \caption{Sampling and Training Times (from \citep{E2ELR}).}
        \label{tab:exp:ED-R:train_time}
    \end{table}

Table \ref{tab:exp:ED-R:inference_time} report the average solving
time using Gurobi (GRB) and the average inference times of the ML
methods.  The ML inference times are reported for a batch of 256
instances on a GPU. The number of gradient steps used by DC3 to
recover feasibility is set to 200 for inference (compared to 50 for
training). For systems with more than 6,000 buses, DC3 is typically
10--30 times slower than other baselines, due to its unrolled gradient
steps.  In contrast, DNN, DeepOPF, and E2ELR all require about 5--10
milliseconds to evaluate a batch of 256 instances.  For the largest
systems, this represents about 25,000 instances per second on a single
GPU. Solving the same volume of instances would take Gurobi more than
a day on a single CPU.

    \begin{table}[!t]
        \centering
            \input{time_inference_DCOPF-R}
        \caption{Solving and Inference Times.}
        \label{tab:exp:ED-R:inference_time}
    \end{table}

\subsection{The Impact of Optimization Proxies}

\begin{figure}[!t]
\centering \includegraphics[width=.75\columnwidth]{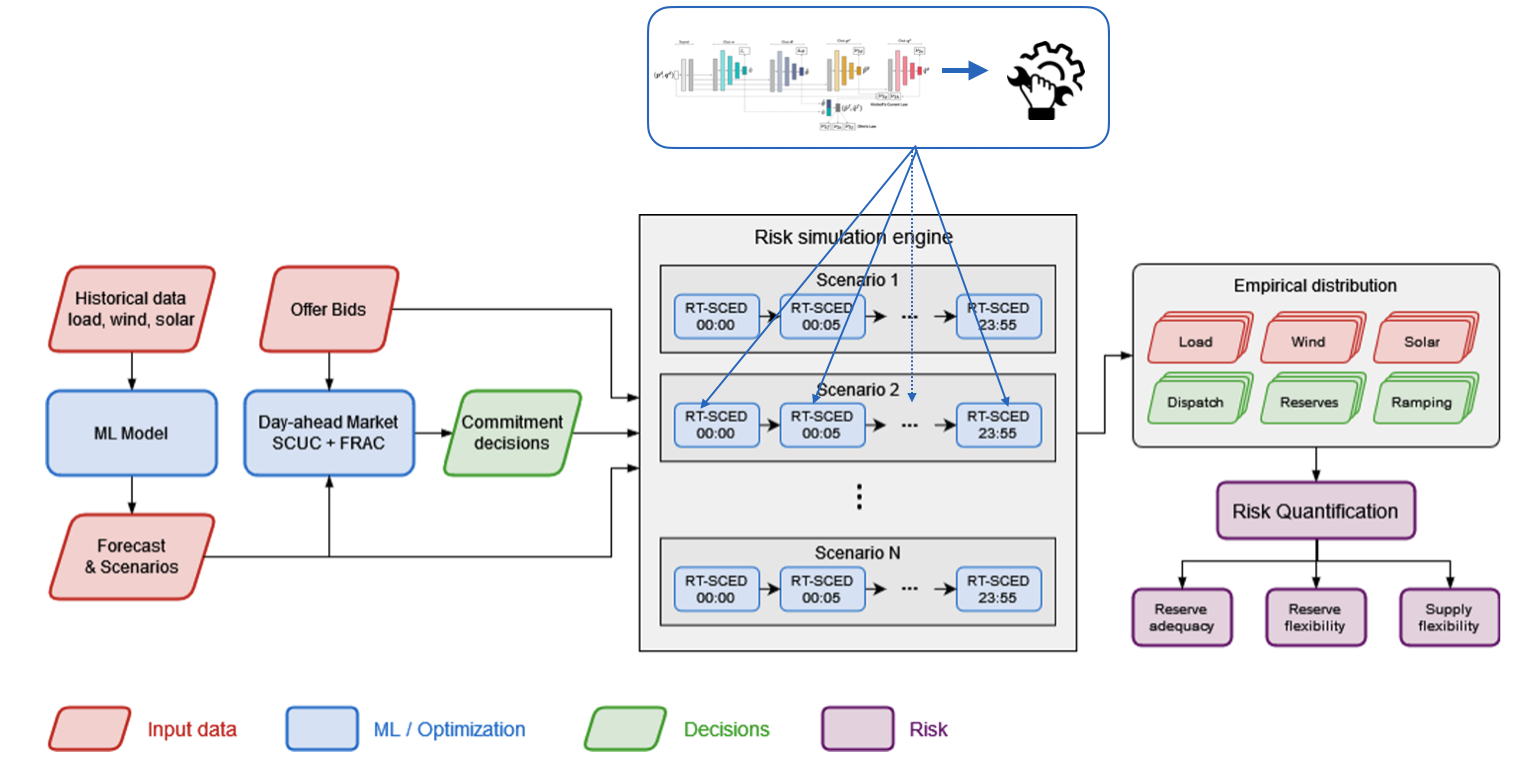}
\caption{A Risk Assessment Framework with Optimization Proxies (from
  \citep{chen2023realtime}).}
\label{fig:ra}
\end{figure}

Optimization proxies have the potential to transform applications as
they bring orders of magnitude improvements in efficiency. This
section highlights these potential benefits on an important topic:
real-time risk assessment. Given the significant growth in volatility
coming from renewable sources of energy, the electrification of the
transportation infrastructure, and data centers, operators are
increasingly worried about the operational and financial risks in
their grids. They are interested in dashboards that give them the
ability to quantify these risks in real time.  Consider, for instance,
the real-time risk-assessment framework described in Figure
\ref{fig:ra}, which runs a collection of Monte-Carlo scenarios for the
next 24 hours. These scenarios are obtained using forecasting methods
(e.g., temporal fusion transformers). The evaluation of each scenario
requires 288 optimizations, i.e., one ED every 5 minutes for the next
24 hours. Solving these 288 optimizations takes about 15 minutes, and a
real-time risk assessment dashboard would need to evaluate thousands
of them.

\begin{figure}[!t]
\centering \includegraphics[width=.75\columnwidth]{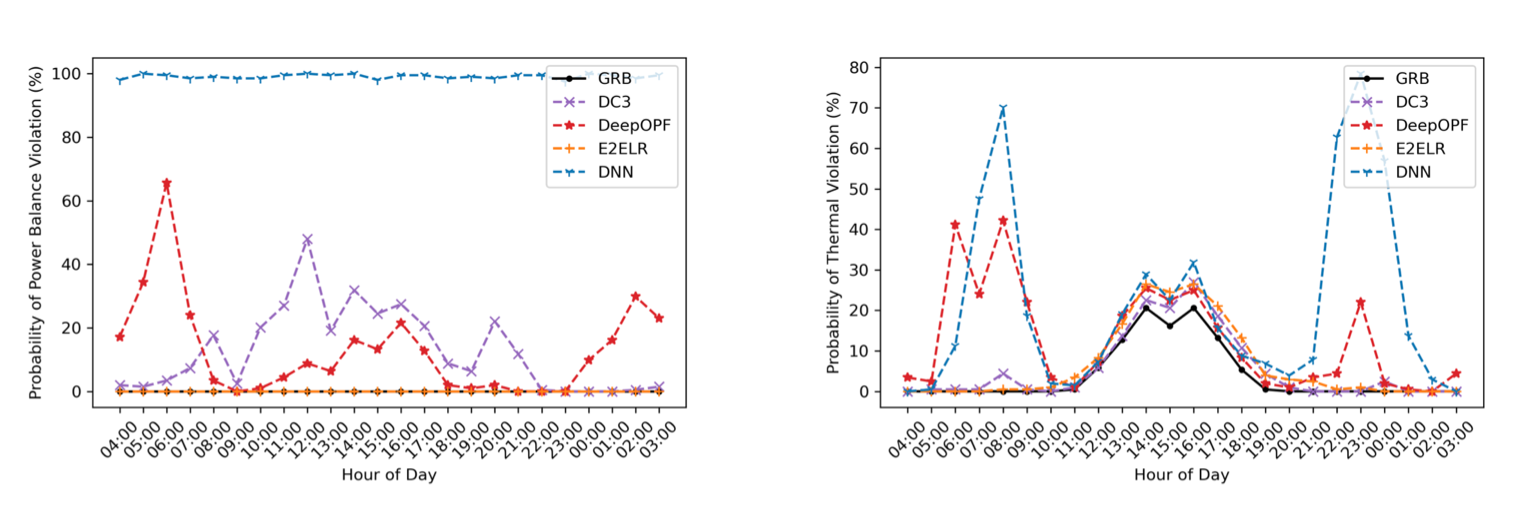}
\caption{Probability of an Adverse Event with a Proxy-based Risk Assessment (from \citep{chen2023realtime}).}
\label{fig:risk}
\end{figure}

{\em Optimization proxies are the key technology enabler to make such
a real-time risk assessment platform a reality.} When primal
optimization proxies replace the optimizations, every scenario is
evaluated in about 5 seconds \citep{E2ELR}, making it possible to
quantify asset-level, system-level, and financial risks in real time
using GPUs. For instance, Figure \ref{fig:risk} reports the probability
of an adverse event for the next 24 hours using E2ELR. It also
compares {\sc E2ELR} with other learning architectures (e.g., DeepOPF,
DC3) and the ground truth, i.e., the risk assessment using a
state-of-the-art optimization solver. The left figure depicts the
probability of an adverse event on the balance constraint: both the
ground truth and {\sc E2ELR} show no such event. In contrast, the
other methods systematically report potential adverse events, since
they cannot guarantee feasibility of that constraint. The right figure
shows the violations of the thermal limit of a congested line. Again,
{\sc E2ELR} isolates the potential violations with great
accuracy. Other methods report many false positives throughout the day.

%% file: gaps_all.tex
\begin{tabular}{llrrrrrrrrr}
\toprule
&
    & \multicolumn{5}{c}{ED-NR}
    & \multicolumn{4}{c}{ED}\\
\cmidrule(lr){3-7}
\cmidrule(lr){8-11}
Loss & System         
        &      DNN &    E2ELR &      DeepOPF & DC3 &     LOOP
        &      DNN &    E2ELR &      DeepOPF & DC3 \\
\midrule
SL   & \ieeeSmall     
        &    69.55 &     \textbf{1.42} &           2.81 &     3.03 &    38.93
        &    75.06 &     \textbf{1.52} &           2.80 &     2.94\\
     & \pegaseSmall   
        &    48.77 &     \textbf{0.74} &           7.78 &     2.80 &    32.53
        &    47.84 &     \textbf{0.74} &           7.50 &     2.97\\
     & \rteLarge           
        &    55.13 &     \textbf{1.35} &          28.23 &     3.68 &    50.21
        &    70.57 &     \textbf{1.82} &          46.66 &     3.49\\
     & \pegaseLarge        
        &    76.06 &     \textbf{0.38} &          33.20 &     1.25 &    33.78
        &    81.19 &     \textbf{0.38} &          30.84 &     1.29\\
     & \pegaseXLarge       
        &    71.14 &     \textbf{0.29} &          64.93 &     1.79 &    52.94
        &    76.32 &     \textbf{0.28} &          69.23 &     1.81\\
     & \gocXXL             
        &   194.13 &     \textbf{0.46} &          55.91 &     2.75 &    36.49
        &   136.25 &     \textbf{0.45} &          41.34 &     2.35\\
\midrule
SSL  & \ieeeSmall          
        &    35.66 &     \textbf{0.74} &           2.23 &     2.51 &    37.78
        &    45.56 &     \textbf{0.78} &           2.82 &     2.80\\
     & \pegaseSmall        
        &    62.07 &     \textbf{0.63} &          10.83 &     2.57 &    32.20
        &    64.69 &     \textbf{0.68} &           9.83 &     2.61\\
     & \rteLarge           
        &    40.73 &     \textbf{1.30} &          42.28 &     2.82 &    50.20
        &    55.16 &     \textbf{1.68} &          48.57 &     3.04\\
     & \pegaseLarge        
        &    43.68 &     \textbf{0.32} &          34.33 &     0.82 &    33.76
        &    44.74 &     \textbf{0.29} &          42.06 &     0.93\\
     & \pegaseXLarge       
        &    57.58 &     \textbf{0.21} &          60.12 &     0.84 &    52.93
        &    61.28 &     \textbf{0.19} &          65.38 &     0.91\\
     & \gocXXL             
        &   108.91 &     \textbf{0.39} &           8.39 &     0.72 &    36.73
        &    93.91 &     \textbf{0.33} &           9.47 &     0.71\\     
\bottomrule
\end{tabular}

%% file: viol_powerBalance_all.tex
\begin{tabular}{llrrrrrr}
\toprule
&
     & \multicolumn{2}{c}{DNN}
     & \multicolumn{2}{c}{DeepOPF}
     & \multicolumn{2}{c}{DC3}\\
\cmidrule(lr){3-4}
\cmidrule(lr){5-6}
\cmidrule(lr){7-8}
Loss & System        
    & \multicolumn{1}{c}{\%feas} & \multicolumn{1}{c}{viol}
    & \multicolumn{1}{c}{\%feas} & \multicolumn{1}{c}{viol}
    & \multicolumn{1}{c}{\%feas} & \multicolumn{1}{c}{viol}\\
\midrule
SL   & \ieeeSmall    
         &      0\% &   0.70  &    100\% &   0.06  &    100\% &   0.25 \\
     & \pegaseSmall  
         &      0\% &   1.90  &     70\% &   1.19  &     65\% &   0.08 \\
     & \rteLarge     
         &      0\% &   2.26  &      1\% &   2.81  &      7\% &   0.01 \\
     & \pegaseLarge  
         &      0\% &   5.91  &      3\% &   5.99  &     29\% &   0.00 \\
     & \pegaseXLarge 
         &      0\% &   7.79  &      1\% &   7.46  &     22\% &   0.00 \\
     & \gocXXL       
         &      0\% &   2.31  &     54\% &   1.54  &     69\% &   0.00 \\
\midrule
SSL  & \ieeeSmall    
         &      0\% &   0.73  &    100\% &   0.56  &    100\% &   0.29 \\
     & \pegaseSmall  
         &      0\% &   2.30  &     63\% &   1.81  &     39\% &   0.05 \\
     & \rteLarge     
         &      0\% &   2.72  &      1\% &   2.51  &      5\% &   0.01 \\
     & \pegaseLarge  
         &      0\% &   7.00  &      3\% &   7.87  &     25\% &   0.08 \\
     & \pegaseXLarge 
         &      0\% &   6.82  &      1\% &   7.83  &     20\% &   0.01 \\
     & \gocXXL       
         &      0\% &   2.61  &     49\% &   1.20  &     62\% &   0.02 \\
\bottomrule
\end{tabular}

%% file: time_inference_DCOPF-R.tex
\begin{tabular}{llrrrrr}
\toprule
Loss & System         &      DNN &    E2ELR  &   DeepOPF &     DC3 & GRB\\
\midrule
SL 
    & \ieeeSmall 
        & \qty[mode=text]{   3.9}{\ms} 
        & \qty[mode=text]{   6.5}{\ms}
        & \qty[mode=text]{4.6}{\ms}
        & \qty[mode=text]{  16.5}{\ms} 
        & \qty[mode=text]{  12.6}{\ms} \\
    & \pegaseSmall 
        & \qty[mode=text]{   4.5}{\ms} 
        & \qty[mode=text]{   6.0}{\ms} 
        & \qty[mode=text]{4.8}{\ms}
        & \qty[mode=text]{  18.9}{\ms} 
        & \qty[mode=text]{  56.5}{\ms} \\
    & \rteLarge 
        & \qty[mode=text]{   5.7}{\ms} 
        & \qty[mode=text]{  10.4}{\ms} 
        & \qty[mode=text]{6.2}{\ms}
        & \qty[mode=text]{  36.1}{\ms} 
        & \qty[mode=text]{ 333.6}{\ms} \\
    & \pegaseLarge 
        & \qty[mode=text]{   6.3}{\ms} 
        & \qty[mode=text]{   7.7}{\ms} 
        & \qty[mode=text]{6.7}{\ms}
        & \qty[mode=text]{  91.6}{\ms} 
        & \qty[mode=text]{1008.0}{\ms} \\
    & \pegaseXLarge 
        & \qty[mode=text]{   8.3}{\ms} 
        & \qty[mode=text]{  10.7}{\ms} 
        & \qty[mode=text]{8.8}{\ms}
        & \qty[mode=text]{ 531.2}{\ms} 
        & \qty[mode=text]{1632.7}{\ms} \\
    & \gocXXL 
        & \qty[mode=text]{   9.3}{\ms} 
        & \qty[mode=text]{  11.1}{\ms}
        & \qty[mode=text]{10.6}{\ms}
        & \qty[mode=text]{ 438.7}{\ms} 
        & \qty[mode=text]{4745.7}{\ms} \\
\midrule
SSL 
    & \ieeeSmall 
        & \qty[mode=text]{   3.9}{\ms}
        & \qty[mode=text]{   7.6}{\ms}
        & \qty[mode=text]{4.4}{\ms}
        & \qty[mode=text]{  17.6}{\ms}
        & \qty[mode=text]{  12.6}{\ms} \\
    & \pegaseSmall 
        & \qty[mode=text]{   4.4}{\ms} 
        & \qty[mode=text]{   5.9}{\ms} 
        & \qty[mode=text]{4.7}{\ms}
        & \qty[mode=text]{  19.1}{\ms} 
        & \qty[mode=text]{  56.5}{\ms} \\
    & \rteLarge 
        & \qty[mode=text]{   6.4}{\ms} 
        & \qty[mode=text]{  10.5}{\ms} 
        & \qty[mode=text]{6.7}{\ms}
        & \qty[mode=text]{  37.3}{\ms} 
        & \qty[mode=text]{ 333.6}{\ms} \\
    & \pegaseLarge 
        & \qty[mode=text]{   7.1}{\ms} 
        & \qty[mode=text]{   8.3}{\ms} 
        & \qty[mode=text]{7.2}{\ms}
        & \qty[mode=text]{  92.9}{\ms} 
        & \qty[mode=text]{1008.0}{\ms} \\
    & \pegaseXLarge 
        & \qty[mode=text]{   7.8}{\ms} 
        & \qty[mode=text]{   8.9}{\ms} 
        & \qty[mode=text]{7.9}{\ms}
        & \qty[mode=text]{ 522.4}{\ms} 
        & \qty[mode=text]{1632.7}{\ms} \\
    & \gocXXL 
        & \qty[mode=text]{  10.2}{\ms} 
        & \qty[mode=text]{  12.4}{\ms} 
        & \qty[mode=text]{10.3}{\ms}
        & \qty[mode=text]{ 435.8}{\ms} 
        & \qty[mode=text]{4745.7}{\ms} \\
\bottomrule
\end{tabular}

%% file: dual_proxies.tex
\section{Dual Optimization Proxies}
\label{section:dual_proxies}

Primal optimization proxies produce (high-quality) feasible solutions
in milliseconds. However, an ideal outcome in optimization practice is
to produce a pair of primal and dual solutions with a small duality
gap. The counterpart for optimization learning would be to produce a
pair of primal and dual proxies that would produce both a feasible
solution and an assessment of its quality in milliseconds. This section
describes a methodology to obtain dual optimization proxies for many
practical applications. To convey the intuition underlying the
methodology, consider the parametric optimization program and its dual

\vspace{-1.2cm}
\begin{multicols}{2}
\begin{flalign*}
\min \qquad & \bc_\x^T \, \y \\
s.t. \qquad & \bA_\x \, \y  = \bb_\x \\
            & \bl_\x \leq \y \leq \bu_x
\end{flalign*}
\break \vspace{-0.5cm}
\begin{flalign*}
\max \qquad & \bb_\x^T \, \z + \bl^T_x \, \z_l - \bu^T_\x \, \z_u \\
s.t. \qquad & \z \, \bA_\x + \z_l - \z_u  = \bc_\x \\
            & \z_l, \z_u \geq 0
\end{flalign*}
\end{multicols}
\noindent
The key insight comes from the fact that the bound constraints $
\bl_\x \leq \y \leq \bu_x $ in the primal model gives a simple way to
find a dual feasible solution: obtain a prediction $\widecheck{\z}$ for
variables $\z$ and use the free dual variables $\z_l$ and $\z_u$ to
obtain a dual feasible solution $(\widecheck{\z},\widecheck{\z}_l,\widecheck{\z}_u)$.
These primal bound constraints are natural in most engineering
applications: for instance, generators in power systems have physical
limits, and the number of trailers in a supply chain is finite.

\begin{figure}[!t]
\centering \includegraphics[width=.75\columnwidth]{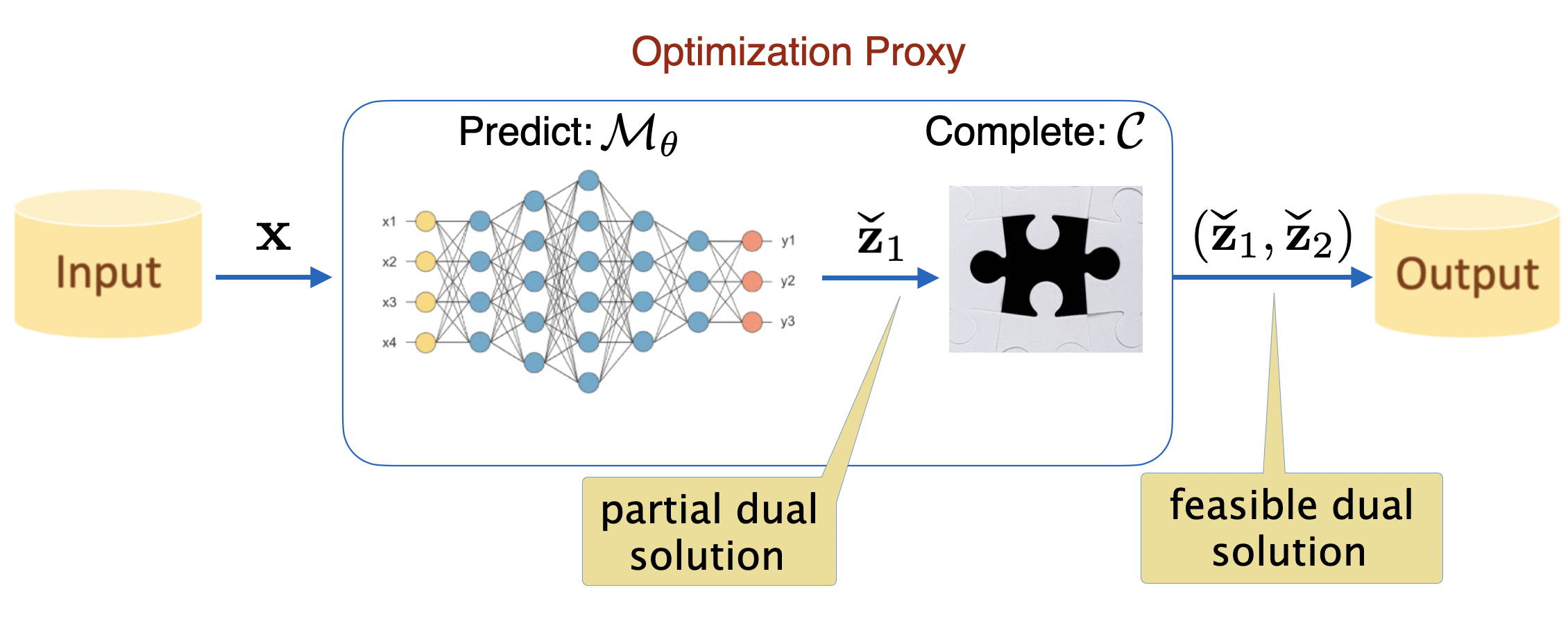}
\caption{The Structure of Dual Optimization Proxies.}
\label{fig:dual_proxy}
\end{figure}

Figure \ref{fig:dual_proxy} provides the general structure of dual
optimization proxies. Given decision variables $\z=(\z_1,\z_2)$, a
dual proxy consists of a parametric machine learning model ${\cal
  M}_{\theta}$ to obtain a prediction $\z_1 = \widecheck{\z}_1$ for a subset of the
decision variables, and a {\em completion} layer ${\cal C}$ to assign
the remaining variables $\z_2 = \widecheck{\z}_2$ and deliver a feasible dual
solution $\widecheck{\z}=(\widecheck{\z}_1,\widecheck{\z}_2)$. The composition
${\cal C} \circ {\cal M}_{\theta}$ is a parametric dual optimization proxy.

Consider now a dual optimization proxy for parametric linear
programming. Given a completion layer ${\cal C}^{lp}$, the composition
${\cal C}^{lp} \circ {\cal M}_{\theta}$ is a parametric dual
optimization proxy that can be trained by the following optimization:
\[
\argmin_{\theta} \frac{1}{n} \sum_{i=1}^n \obji((\mathcal{C}^{lp} \circ {\cal M}_{\theta})(\x_i)).
\]
where $\obji((\z,\z_l,\z_u)) = \bb_{\x_i}^T \, \z + \bl_{\x_i}^T \,
\z_l - \bu^T_{\x_i} \, \z_u$.

It remains to derive the completion layer ${\cal C}^{lp}$. Given a prediction $\widecheck{\z}^t
= {\cal M}_{\theta^t}(\x)$ at step $t$, the dual problem becomes
\begin{flalign*}
\max \qquad & \bl^T_\x \, \z_l^t - \bu^T_\x \, \z_u^t \\
s.t. \qquad &  \z_l^t - \z_u^t  = \bc_\x - \widecheck{\z}^t \, \bA_\x \\
            & \z_l^t, \z_u^t \geq 0
\end{flalign*}
\noindent
The optimal solution to this problem has a closed form \citep{klamkin2024dual}:
\begin{subequations}
\label{eq:maxzlzu}  
\begin{align*}
       \widecheck{\bz}_l^t &= |\bc_\x - \widecheck{\z}^t \, \bA_\x|^{+} \\
       \widecheck{\bz}_u^t &= |\bc_\x - \widecheck{\z}^t \, \bA_\x|^{-} 
\end{align*}
\end{subequations}
where $|v|^{+} \, {=} \, \max(0, v)$, $|v|^{-} \, {=} \, \max(0, -v)$,
and $(\widecheck{\z}^t,\, \widecheck{\bz}_l^t,\, \widecheck{\bz}_u^t)$ is a
feasible solution to the dual problem. Hence, the completion layer
${\cal C}^{lp}$ for parametric linear programs is defined as
\[
  {\cal C}^{lp}(\widecheck{z}) = (\widecheck{z},|\bc_\x - \widecheck{\z} \, \bA_\x|^{+},|\bc_\x - \widecheck{\z} \, \bA_\x|^{-}).
  \]
The forward pass of the dual proxy ${\cal C}^{lp} \circ {\cal
  M}_{\theta}$ can be evaluated efficiency. Moreover, the dual proxy
admits subgradients everywhere, enabling a fast backward pass as well.

Dual optimization proxies have been proposed first for the
Second-Order Cone Programming (SOCP) relaxation of the AC power flow
equations \citep{DualProxies}.  The completion layer is
particularly interesting in that paper: it uses properties of the dual
optimal solutions to determine how to complete the set of dual
variables not predicted by the neural network. Again, the dual
optimization proxy is a differentiable program that can be trained
end-to-end to produce feasible solutions at training and inference
times. Experimental results have shown that the resulting proxies can
find near-optimal dual feasible solutions. Dual optimization proxies
for conic optimization have been studied extensively in
\citep{tanneau2024dual}. More generally, a methodology to find quality
guarantees for an optimization problem $\Phi$ can be summarized as
follows:
\begin{enumerate}[noitemsep,topsep=2.5pt,parsep=0pt,partopsep=0pt,leftmargin=0.5cm]
\item find a conic approximation $\Psi$ to problem $\Phi$;
\item obtain the dual $\Psi_d$ of $\Psi$;
\item derive a dual optimization proxy ${\Psi}_d^\downarrow$ for $\Psi_d$.
\end{enumerate}

\subsection{Experimental Results}

This section presents some experimental result on the self-supervised
Dual Optimal Proxy for Linear Programming (DOPLP) applied to the
DC-OPF problem. The full details and results can be found in
\citep{klamkin2024dual}. The DC-OPF formulation is shown in Model
\ref{model:DC-OPF}
\begin{model}[!t]
\begin{subequations}
    \label{eq:DCOPF}
    \begin{flalign}
        \argmin_{\pg,\pf} \quad & \mathbf{c}^{\top} \pg & & \\
        s.t. \quad 
        & \mathbf{e}^{\top} \p = \mathbf{e}^{\top} \pd & & \\
        & \pf = \mathbf{PTDF}(\pg - \pd) & & \\
        & \pgmin \leq \pg \leq \pgmax & & \\
        & \pfmin \leq \pf \leq \pfmax & &
    \end{flalign}
\end{subequations}
\caption{The DC-OPF Model.}
\label{model:DC-OPF}
\end{model}
It is close to the ED: it does not have reserve constraints but its
line constraints are hard.

DOPLP uses a 3-layer fully connected neural network with ReLU
activations. The training uses 10,000 feasible instances and 2,500
instances for validation. The testing set is the same set of 5,000
feasible instances. The experiments are repeated 10 times with
different seeds to ensure results do not depend on a particular
training/validation split or initial neural network weights.  All
training was run on NVIDIA Tesla V100 16GB GPUs.

\begin{table}[!t]
    \centering
    \vskip 0.1in
    \begin{tabular}{l|cc|cc}
        \toprule
        Benchmark & $m$ & $n$ & $h$ & $|\theta|$\\
        \midrule
        \texttt{1354\_pegase} & 1992 & 2251 & 2048 & 16.6M\\
        \texttt{2869\_pegase} & 4583 & 5092 & 4096 & 71.1M\\
        \texttt{6470\_rte} & 9006 & 9766 & 8192 & 281M\\
        \bottomrule
    \end{tabular}
    \caption{DCOPF Benchmark Summary}
    \label{tab:dcopf_benchmark}
    \vskip -0.1in
\end{table}

DOPLP is evaluated on three large-scale DC Optimal Power Flows (DCOPF)
instances from the PGLib library \citep{PGLib}: \texttt{1354\_pegase},
\texttt{2869\_pegase}, and \texttt{6470\_rte}. The instances are
generated as in Section \ref{section:experiments}, except for
\texttt{6470\_rte} where $U$ is set to 1.05 instead, since the given
instance is already congested. Table \ref{tab:dcopf_benchmark} reports
the number of constraints $m$, the number of variables $n$, the hidden
layer size $h$, and the total number of parameters $|\theta|$ in the
neural network ${\cal M}_\theta$ for each benchmark. For the purposes
of assessing the quality of the results, each instance is solved using
Mosek 10.1.24 \citep{mosek} to obtain the optimal dual solution
$\by^*$.

Table \ref{tab:dual_gap} shows the minimum, geometric mean,
$99^\text{th}$ percentile, and maximum dual gap ratio over the testing
set samples.  Each entry in the table is the geometric mean $\pm$ the
standard deviation over the ten trials. {\em DOPLP produces mean
optimality gaps in the range of $[0.25,0.5]\%$, showing that DOPLP is
capable of producing dual feasible solutions with extremely tight dual
bounds for large-scale DCOPF problems.}

\begin{table*}[!t]
    \centering
    \vskip 0.1in
    \begin{tabular}{l|cccccc}
        \toprule
        Benchmark & Minimum & Geometric Mean & $99^\text{th}$ Percentile & Maximum \\
        \midrule
        \multirow{1}{*}{\texttt{1354\_pegase}} & 
          $ 0.14\%{\scriptstyle (\pm .014\% )}$ & $ 0.23\%{\scriptstyle (\pm .016\% )}$ & $ 2.15\%{\scriptstyle (\pm .030\% )}$ & $ 13.76\%{\scriptstyle (\pm .015\% )}$ \\
        \midrule
        \multirow{1}{*}{\texttt{2869\_pegase}} &
          $ 0.15\%{\scriptstyle (\pm .018\% )}$ & $ 0.25\%{\scriptstyle (\pm .106\% )}$ & $ 3.73\%{\scriptstyle (\pm .194\% )}$ & $ 10.26\%{\scriptstyle (\pm .513\% )}$ \\
        \midrule
         \multirow{1}{*}{\texttt{6470\_rte}} &
          $ 0.15\%{\scriptstyle (\pm .013\% )}$ & $ 0.50\%{\scriptstyle (\pm .146\% )}$ & $ 1.96\%{\scriptstyle (\pm .615\% )}$ & $ 6.95\%{\scriptstyle (\pm .633\% )}$ \\
        \bottomrule
    \end{tabular}
    \caption{Dual Gap Ratio (adapted from \citep{klamkin2024dual}).}
    \label{tab:dual_gap}
\end{table*}

%% file: primal_dual_proxies.tex
\section{Primal-Dual Learning} 
\label{section:PDL}

Another approach to find feasible solutions is to develop optimization
proxies that adapt traditional optimization algorithms to the learning
context.  Consider the constrained optimization problem
\begin{equation*}
\label{eq:opt-PDL}  
\argmin_{\y} \obj(\y) \mbox{ subject to }  \eq(\y) = 0. 
\end{equation*}
\noindent
The {\em Augmented Lagrangian Method} (ALM) solves unconstrained
optimization problems of the form
\begin{equation}
\label{eq:ALM}
\argmin_{\y} \obj(\y) \!+\! \deq^T \, \eq(\y) + \frac{\rho}{2} \, \mathbf{1}^{\top} \eq(\y)^2
\end{equation}
where $\rho$ is a penalty coefficient and $\deq$ are the Lagrangian
multiplier approximations. These multipliers are updated using the
rule
\begin{equation}
\label{eq:UR}
\deq \leftarrow \deq + \rho \, \eq(\y).
\end{equation}

\begin{figure}[!t]
\centering \includegraphics[width=.80\columnwidth]{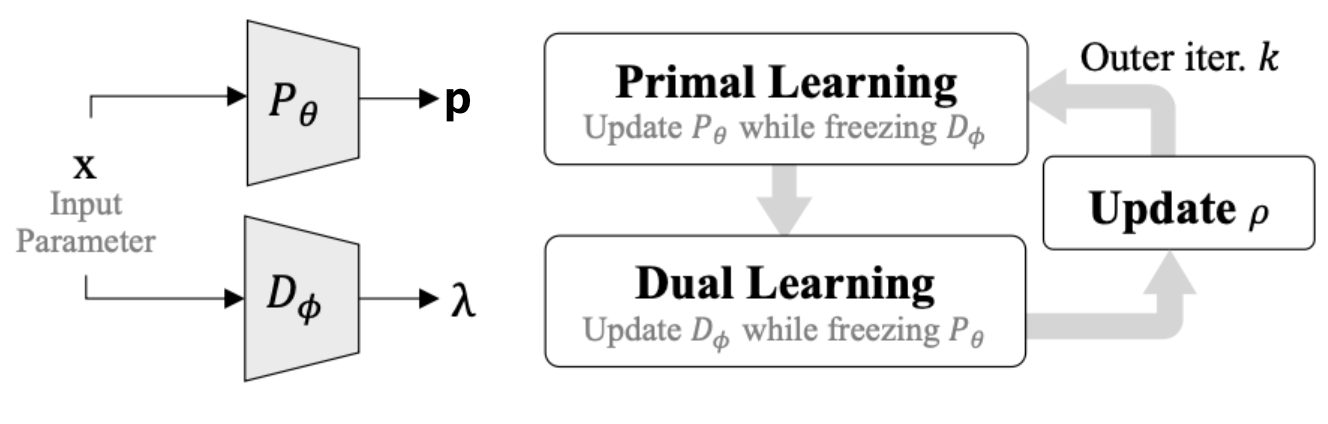}
\caption{The Architecture of Primal-Dual Learning (from \citep{PDL}).}
\label{fig:PDL}
\end{figure}

\noindent
{\em Primal-Dual Learning} (PDL) \citep{PDL} is a self-supervised
method that jointly learns two neural networks: a primal neural
network $\pnet$ that learns the input/output mapping of the ALM
unconstrained optimizations and a dual network $\dnet$ that learns the
dual updates. The overall architecture of PDL is shown in Figure
\ref{fig:PDL}. At each iteration, the {\em primal learning} step
updates the learning parameters $\theta$ of the primal network while
keeping the dual network $\dnet$ fixed. After completion of the primal
learning, PDL applies a {\em dual learning} step that updates the
learnable parameters $\phi$ of the dual network $\dnet$.

More precisely, at step $t$, PDL has the values $\theta^t$ and $\phi^t$ for the primal
and dual learning parameters and seeks to improve them into $\theta^{t+1}$ and $\phi^{t+1}$.
The primal learning problem of step $t$ is the optimization
\begin{equation}
  \label{eq:PDL_primal}
\theta^{t+1} = \argmin_{\theta} \frac{1}{n} \sum_{i \in [n]} \obji(\y_i) + \!+\! D_{\phi^t}(\x_i)^T \, \eqi(\y_i) + \frac{\rho}{2} \, \mathbf{1}^{\top} \eqi(\y_i)^2,
\end{equation}
where $\y_i = P_{\theta}(\x_i)$. The dual learning problem is the optimization
\begin{equation}
\label{eq:PDL_dual}
\phi^{t+1} = \argmin_{\phi} \frac{1}{n} \sum_{i \in [n]} \norm{D_{\phi}(\x_i) \!-\! D_{\phi^t}(\x_i) \!+\! \rho \, \eqi(P_{\theta^{t+1}}(\x_i))}.
\end{equation}

\noindent
These two steps, i.e., the training of the primal and dual networks
that mimic the ALM, are iterated in sequence until convergence. To
handle severe violations, each iteration may increase the penalty
coefficient $\rho$.  At iteration $t$, this update uses the maximum
violation $v^{t+1}$
\begin{equation}
\label{eq:rho_sigma}
v^{t+1} = \max_{i \in [n]}\left\{\linfnorm{\eq(P_{\theta^{t+1}}(\x_i))}\right\}.
\end{equation}
The penalty coefficient increases when the maximum violation $v^{t+1}$ is
greater than a tolerance value $\tau$ times the maximum violation from
the previous iteration $v^{t}$, i.e.,
\begin{equation}
\label{eq:update_rho-main}
\rho \leftarrow \min\{\alpha\rho,\rhomax\}\;\text{if } v^{t+1} > \tau v^{t},
\end{equation}
where $\tau\in(0,1)$ is the tolerance, $\alpha\!>\!1$ is an update
multiplier, and $\rhomax$ is an upper bound on the penalty
coefficient. The PDL algorithm is summarized in Algorithm \ref{alg:pdl}.

\begin{algorithm}[!t]
\textbf{Parameter}: 
    Initial penalty coefficient $\rho$, 
    Maximum outer iteration $T$, 
    Penalty coefficient updating multiplier $\alpha$, 
    Violation tolerance $\tau$, 
    Upper penalty coefficient safeguard $\rhomax$ \\ 
\textbf{Input}: Input Distribution ${\cal P}$\\
\textbf{Output}: learned primal and dual nets $\pnet$, $\dnet$
\begin{algorithmic}[1] 
\State Initialize $\theta^1$ and $\phi^1$ randomly;
\For{$t \in [T]$}
    \State Train the primal network using Equation \eqref{eq:PDL_primal} 
    \State Train the dual network using Equation \eqref{eq:PDL_dual}
    \State Update $\rho$ using Equation \eqref{eq:update_rho-main}
\EndFor
\State \textbf{return} $P_\theta = P_{\theta^{T+1}}$ and $D_{\phi} = D_{\phi^{T+1}}$
\end{algorithmic}
\caption{The Primal-Dual Learning Algorithm (PDL).}
\label{alg:pdl}
\end{algorithm}

There is a fundamental difference between PDL and the Lagrangian dual
method presented in Section \ref{section:erm}. In these Lagrangian
methods, the constraint multipliers used to balance feasibility and
optimality are not instance-specific: they are aggregated for each
constraint over all the instances. This aggregation limits their
capability to ensure feasibility. In contrast, by jointly learning the
primal and dual networks, PDL is capable of using {\em
  instance-specific multipliers} for penalizing constraints. As a
result, during each iteration, the training process learns the primal
and dual iteration points of the ALM algorithm. Eventually, these
iteration points, and hence the primal and dual networks, are expected
to converge to the primal and dual solutions of the optimization
problem.

{\em Primal-Dual Optimization proxies have the ability to deploy
  optimization models that would otherwise be too complex to meet
  real-time requirements.} By shifting most of the computational
burden offline during training, optimization proxies provide a way to
produce high-quality solutions to optimization models that cannot be
solved fast enough to meet real-time constraints. An interesting
example is the Security-Constrained Optimal Power Flow problem (SCOPF)
that captures the automatic primary response of generators in case of
transmission line or generator contingencies. The next subsection
briefly describes this application in some detail. The presentation
follows \citep{PDLSCOPF}, where the full results and technical
description can be found.  An interesting feature of this Primal-Dual
Learning application is that it mimics the Column and Constraint
Generation Algorithm (CCGA) in
\citep{velloso2021exact,velloso2021combining}, which is the state of
the art for the SCOPF.

\subsection{The SCOPF Problem}
\label{ssec:SCOPF}

\begin{model}[!t]
\caption{The Extensive SCOPF Formulation.}
\label{model:scopf_ext}
\begin{flalign}
    \:\:\:\:\:\:\:\:
    \:\:\:\:\:
    \mmin_{
    \mathclap{
    \substack{ \setlength{\jot}{-0.8\baselineskip}\everymath{\scriptstyle}
    \begin{array}{c}
    \g,\\ 
    \!\!\left[\gk\!,\rhok\!,\nk\!\right]_{k\in\Kg}\!, \\
    \!\!\left[\slack_k\right]_{\!k\in\{\!0\!\}\cup\Kg\cup\Ke}
    \end{array}
  }
  }
  }
 \:\:\:\:& \objc^{\top}\g \!+\! \slackpen\!\!\left(\sum_{k\in\{\!0\!\}\cup\Kg\cup\Ke}\!\!\norm{\slack_k}_{\!1} \!\!\right)\label{eq:scopf_ext_obj}\\   
\text{\textbf{s. t.}:}\:\:\:\:\:& \mathbf{1}^\top\g = \mathbf{1}^\top\load \label{eq:scopf_ext_cnst_pb} \\
 & \flb\!-\!\slack_0 \leq \f\!=\! \mathbf{PTDF}(\load\!-\!\g) \leq \fub\!+\!\slack_0  \label{eq:scopf_ext_cnst_pf} \\
 & \glb \leq \g \leq \gub \label{eq:scopf_ext_cnst_gbound} \\
 & \mathbf{e}^\top\gk = \mathbf{e}^\top\load & & \forall k\!\in\!\Kg \label{eq:scopf_ext_cnst_pb_kg} \\
& \flb\!-\!\slack_k \!\leq\! \f_k\!=\!\mathbf{PTDF}(\load\!-\!\gk) \!\leq\! \fub\!+\!\slack_k& &\forall k\!\in\!\Kg \label{eq:scopf_ext_cnst_pf_kg} \\
 & \underline{p}_i \leq p_{k\!,i} \leq \overline{p}_i                  & & \forall i\!\in\!\cG,\forall k\!\in\!\Kg,i\!\neq\! k \label{eq:scopf_ext_cnst_gbound_kg} \\
 & p_{k\!,k} = 0 & & \forall k\!\in\!\Kg \label{eq:scopf_ext_cnst_gbound_kg2} \\
 & |p_{k\!,i}\!-\!p_i\!-\!n_k\gamma_i\ddot{p}_i|\!\leq\!\ddot{p}_i\rho_{k\!,i} & & \forall i\!\in\!\cG,\forall k\!\in\!\Kg,i\!\neq\! k \label{eq:scopf_ext_cnst_apr1} \\
 & p_i\!+\!n_k\gamma_i\ddot{p}_i\!\geq\!\ddot{p}_i\rho_{k\!,i}\!+\!\underline{p}_i & & \forall i\!\in\!\cG,\forall k\!\in\!\Kg,i\!\neq\! k \label{eq:scopf_ext_cnst_apr2} \\
 & p_{k\!,i} \!\geq\!\ddot{p}_i\rho_{k\!,i}\!+\!\underline{p}_i & &\forall i\!\in\!\cG,\forall k\!\in\!\Kg,i\!\neq\! k \label{eq:scopf_ext_cnst_apr3} \\
 & \flb\!-\!\slack_k \leq \f \!+\! f_k \LODF_k  \leq \fub\!+\!\slack_k & & \forall k\!\in\!\Ke \label{eq:scopf_ext_cnst_pf_ke} \\
 & \slack_k \geq 0 & & \forall k\!\in \{\!0\!\}\!\cup\!\Kg\!\cup\!\Ke \label{eq:scopf_ext_cnst_slack_bound} \\ 
 & \nk \in \left[0,1\right] & &\forall k\!\in\!\Kg \label{eq:scopf_ext_cnst_n_kg_bound} \\
 & \rho_{k\!,i} \in \left\{0,1\right\} & & \forall i\!\in\!\cG,\forall k\!\in\!\Kg,i\!\neq\! k \label{eq:scopf_ext_cnst_rho_kg_bound}
\end{flalign}
\end{model}

Model~\ref{model:scopf_ext} presents the extensive formulation of the
SCOPF with $N-1$ generator and line contingencies. The formulation is
similar in spirit to the ED, but it replaces the reserve constraints
with actual generator and line contingencies. The primary objective of
the SCOPF is to determine the generator setpoints for the base case
while ensuring the feasibility of each generator and line
contingencies.  In other words, even if a generator or line fails, the
base setpoints should provide enough flexibility to be transformed
into a feasible solution for the contingency. The
objective~\eqref{eq:scopf_ext_obj} sums the linear cost of the base
case dispatch $\g$ and the penalties for violating the thermal limits
in the base case and in the contingencies.

The base case includes Constraints~\eqref{eq:scopf_ext_cnst_pb},
\eqref{eq:scopf_ext_cnst_pf}, and \eqref{eq:scopf_ext_cnst_gbound}.
They include the hard power balance constraints in the base case, the
soft thermal limits, and the generation limits.  Each generator
contingency imposes Constraints~\eqref{eq:scopf_ext_cnst_pb_kg},
\eqref{eq:scopf_ext_cnst_pf_kg}, and
\eqref{eq:scopf_ext_cnst_gbound_kg} to enforce the power balance, the
soft thermal limits, and the generation bounds under the generator
contingency.  The only difference from the base case is
Constraint~\eqref{eq:scopf_ext_cnst_gbound_kg2} that specifies that
generator $k$ should remain inactive under its contingency.

It is important to model accurately how the generators react after a
contingency. In Model~\ref{model:scopf_ext}, their response follows an
Automatic Primary Response (APR) control mechanism
\citep{dvorkin2016optimizing,aravena2022recent}, as shown in
Constraints~\eqref{eq:scopf_ext_cnst_apr1}--\eqref{eq:scopf_ext_cnst_apr3}. The
formulation here is based on the APR model used in
\citep{velloso2021exact,velloso2021combining} where, for each
generator contingency $k$, a system-wide signal $n_k \in[0,1]$
represents the level of system response required to resolve the power
imbalance.  The APR control also ensures that the change in dispatch
under a contingency is proportional to the droop slope, which is
determined by the product of the generator capacity $\ddot{p}$ and the
predefined parameter $\gamma$ as in \citep{velloso2021exact}.  The
generator capacities are defined as $\ddot{p}\!=\!\gub\!-\!\glb$.  The
APR constraints ensure that the generation dispatch under generator
contingency remains within the generation limits, leading to the
following generator response:
\begin{equation}
\label{eq:apr}
p_{k,i} = \min\{p_i\!+\!n_k\gamma_i\ddot{p}_i,\overline{p}_i\},\,\,\forall i\!\in\!\cG, \forall k\!\in\!\Kg,i\!\neq\!k.
\end{equation}
The binary variables $\rho_{k,i}$ implement the disjunctive
constraint~\eqref{eq:apr} for all $i\in\cG$ and $k\in\Kg$ such that
$k\neq i$. Line contingencies are handled by
Constraints~\eqref{eq:scopf_ext_cnst_pf_ke}.  These constraints ensure
that, during a line contingency, there is an immediate redistribution
of power flow specified by the Line Outage Distribution Factor (LODF)
\citep{guo2009direct,tejada2017security}.  Under line contingency $k$,
the $k$-th column vector $\LODF_k$ of the LODF matrix
$\LODF$, delineates the redistribution
of base case power flow $f_k$ at line $k$ to the other lines to ensure
that there is no power flow at line $k$. Note that
Model~\ref{model:scopf_ext} is a large-scale MIP model.

\subsection{Primal-Dual Learning for SCOPF}
\label{section:PDL-SCOPF}

\begin{figure}[!t]
\centering
\includegraphics[width=.99\columnwidth]{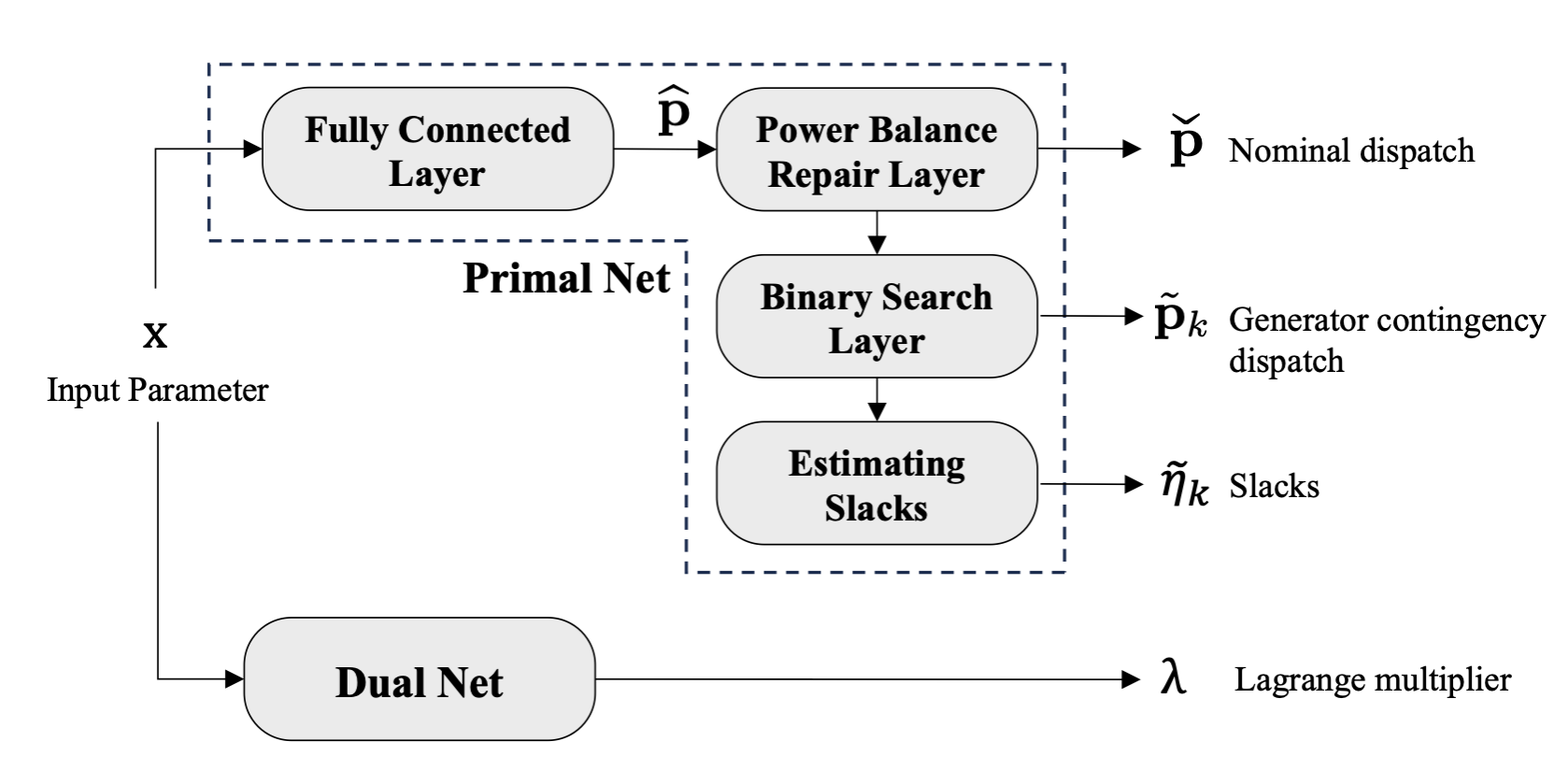}
\caption{The Primal and Dual Networks of PDL-SCOPF (Adapted from \citep{PDLSCOPF}.}
\label{fig:schematic_pdl}
\end{figure}

This section describes PDL-SCOPF, i.e., the Primal-Dual Learning of
large-scale SCOPFs. The primal variables to estimate are $\y:=\p,
\{\gk,n_k,\rho_k\}_{k \in \Kg},
\{\slackapp_k\}_{\{0\}\cup\Kg\cup\Ke}$, the objective function
$\obj(\y)$ is the original objective function \eqref{eq:scopf_ext_obj}
of Model \ref{model:scopf_ext}, and the constraints $\eq(\y)$ capture
the power balance equations \eqref{eq:scopf_ext_cnst_pb_kg} for the
generator contingencies. Figure~\ref{fig:schematic_pdl} provides a schematic
representation of the primal and dual networks which, given the input
configuration vector $\x$, estimate the primal and dual solutions for
SCOPF. The design of the primal learning network of PDL-SCOPF has two
notable features: (1) a repair layer for restoring the power balance
of the base case like in E2ELR, and (2) a binary search layer to
estimate the generator dispatches in the generator contingencies
inspired by the CCGA. 

\subsection{The Primal and Dual Networks}

The primal network estimates the nominal dispatch, the contingency
dispatches, and the slacks of the thermal constraints. As shown in
Figure~\ref{fig:schematic}, it uses three main components: (1) a fully
connected layer with a sigmoid layer that produces a first
approximation $\widehat{\p}$ of the base dispatch; (2) a repair layer that
produces a feasible dispatch $\widecheck{\p}$ to the base dispatch;
and (3) a binary search layer that computes an approximation
$\tilde{\p}_k$ of dispatch for contingency $k$ by mimicking the binary
search of the CCGA. {\em These three components, and the computation
  of the constraint slacks, constitute a differentiable program for
  the primal learning step that is trained end-to-end.}

\begin{algorithm}[!t]
\caption{The Binary Search Layer \texttt{BSLayer}($\g$,$k$) (Adapted from \citep{PDLSCOPF}).}
\label{alg:bslayer}
\textbf{Parameter}:  Maximum iteration $t$\\
Initialize: $n_k\!=\!0.5$, $n_{\text{min}}\!=\!0$, $n_{\text{max}}\!=\!1$
\begin{algorithmic}[1] 
\For{$j=0,1,\dots,t$}
    \State $p^{(j)}_{k\!,i} \leftarrow \min\{p_i\!+\!n_k\gamma_i\ddot{p}_i, \overline{p}_i\}$, $\forall i\in\cG$
    \State $p^{(j)}_{k\!,k} \leftarrow 0$
    \State $e_k\leftarrow \mathbf{1}^{\top}\!\mathbf{p}^{(j)}_k \!-\! \mathbf{1}^{\top}\!\load$
    \State \textbf{if} $e_k>0$ \textbf{then:} $n_{\text{max}}\leftarrow n_k$
    \State \textbf{else:} $n_{\text{min}}\leftarrow n_k$
    \State $n_k \leftarrow 0.5(n_{\text{max}}\!+\!n_{\text{min}})$
\EndFor
\For{$i\in\cG$}
\State \textbf{if} $p_i\!+\!n_k\gamma_i\ddot{p}_i > \overline{p}_i$ \textbf{then:} $\rho_{k\!,i}\leftarrow 1$
\State \textbf{$\!$else:} $\rho_{k\!,i}\leftarrow 0$
\EndFor
\State\textbf{return} $\mathbf{p}^{(j)}_k,\nk,\rhok$
\end{algorithmic}
\end{algorithm}

Estimating the dispatches for all $N\!-\!1$ generator contingencies
presents a computational challenge, as the number of binary variables
grows quadratically with the number of generator contingencies.  For
this reason, the primal network only predicts the base dispatch and
uses a binary search layer, inspired by the CCGA algorithm, to compute
the generator dispatches under all contingencies. {\em In other words,
  the PDL proxy only predicts the base case dispatch and computes the
  contingency dispatches through an efficient implicit layer.} The
{\em binary search layer} is an adaptation of its CCGA counterpart and
is described in Algorithm~\ref{alg:bslayer}. It estimates the
dispatches under the generator contingencies ($\gk$) and the
APR-related variables ($\nk, \rhok$) from the base dispatches. The
algorithm performs a binary search on the global signal $n_k$ in order
to try to find contingency dispatches that satisfy the power balance
constraint. Contrary to its CCGA counterpart,
Algorithm~\ref{alg:bslayer} may not always satisfy the power balance
constraints in the contingencies for reasons described in
\citep{PDLSCOPF}; They are however satisfied when training completes
in the experiments.

The forward pass of the binary search layer is conducted using
Algorithm~\ref{alg:bslayer}. Its backward pass computes the gradient
of the contingency dispatches, which are expressed in terms of 
the base case dispatches and the global signal $n_k$ computed by
the binary search in the forward pass, i.e.,
\[
\tilde{p}_{k,i} = \min\{\widecheck{p}_i\!+\!n_k\gamma_i\ddot{p}_i,\overline{p}_i\},\,\,\forall i\!\in\!\cG, \forall k\!\in\!\Kg,i\!\neq\!k.
\]
This expression obviously has subgradients everywhere.

Once the generation dispatches for the base case and generator
contingencies are estimated, it is possible to calculate the slack
variables for the base case, the generator contingencies
\eqref{eq:scopf_ext_cnst_pf_kg}, and the line contingencies
\eqref{eq:scopf_ext_cnst_pf_ke} easily.

The PDL proxy has the nice property that almost all constraints are
implicitly satisfied. Only the power balance constraints in the
generator contingencies may be violated. It is precisely those
constraints that are captured in the primal loss function of
PDL-SCOPF, i.e.,
\begin{equation*}
\eq(\y)_k = \mathbf{1}^{\top}\!\mathbf{\tilde{p}}_k \!-\! \mathbf{1}^{\top}\!\load \;\;\; (k\in\Kg).
\end{equation*}
As a result, since all the other constraints are satisfied, the dual
network produces optimal dual estimates $\deq=\dnet(\x)$ for the
generator contingency power balance constraints.

\subsection{Numerical Experiments}

The full details of the experiments are given in
\citep{PDLSCOPF}. This section summarizes some critical points.  The
effectiveness of PDL-SCOPF is assessed on five specific cases from the
Power Grid Library (PGLIB) \citep{babaeinejadsarookolaee2019power}
given in Table~\ref{tab:case_spec}.

\begin{table}[!t]
\centering
\small
\setlength{\tabcolsep}{4pt}
\begin{tabular}{@{}lccccccccc@{}}
\toprule
Test Case  & $|\cN|$ & $|\cG|$ & $|\cL|$ & $|\cE|$ & $\lvert\Kg\rvert$ & $\lvert\Ke\rvert$ &$\dim\!\left(\x\right)$\\
\midrule
\texttt{1354\_peg}     & 1354   & 260  & 673  & 1991    & 193   & 1430   & 1193 \\
\texttt{1888\_rte}     & 1888   & 290  & 1000 & 2531    & 290   & 1567   & 1580 \\
\texttt{3022\_goc}     & 3022   & 327  & 1574 & 4135    & 327   & 3180   & 2228 \\
\texttt{4917\_goc}     & 4917   & 567  & 2619 & 6726    & 567   & 5066   & 3753 \\
\texttt{6515\_rte}     & 6515   & 684  & 3673 & 9037    & 657   & 6474   & 5041 \\
\bottomrule
\end{tabular}
\caption{Specifications of the SCOPF Test Cases.}
\label{tab:case_spec}
\end{table}

\paragraph{Instance Generation}

The instance generation perturbs the load demands, the cost
coefficients, and the upper bounds of the generation dispatch, i.e.,
$x:=\{\load,\objc,\gub\}$, in the PGLIB configuration.  This
generalizes earlier settings \citep{AAAI2020,DC3,PDL}, where
the only load demand $\load$ is perturbed. The load demands were
sampled from the input distribution $\cal I$, which is defined as the
truncated multivariate Gaussian distribution.  To perturb the cost
coefficients and the dispatch upper bounds, the experiments use the
base values already provided in PGLIB and multiply them by factors
specific to each instance.  Table~\ref{tab:numvar} highlights the size
of the extensive SCOPF formulation~\eqref{model:scopf_ext}: it gives
the numbers of variables and constraints before and after applying the
Presolve (with the default setting) in Gurobi. Presolve cannot even
complete due to out of memory for the larger test cases.  For the
\texttt{3022\_goc}, there are about 14.6 million continuous variables
and 106,600 binary variables. In the largest case \texttt{6515\_rte},
there are about 64.9 million continuous variables, 449,400 binary
variables, and 130.7 million constraints.

The evaluation uses 1,000 instances sampled from the predefined
distribution. Optimal solutions were obtained using CCGA
\citep{velloso2021exact}. Table~\ref{tab:CCGA_time}
presents the computing times and the number of iterations of the
CCGA. Observe the challenging \texttt{6515\_rte} case, which requires
at least 4648 seconds to solve.

PDL-SCOPF has an interesting feature: the instances used in training
are generated on the fly from the distribution ${\cal P}$ at each
iteration of Algorithm \ref{alg:pdl}. There is no need to generate a
set of instances before training.

\paragraph{Baselines}

PDL-SCOPF is evaluated against three machine learning baselines to
assess its performance.  {\em It is important to stress that all
  baselines use the same primal neural architecture, including the
  repair layer for the power balance in the base case.}  Only
PDL-SCOPF also uses a dual network.  The first baseline, denoted as
\emph{Penalty}, is a self-supervised framework that uses the function
$ \obj(\y) + \rho \mathbf{1}^{\top}\viol{\eq(\y)}$ as a loss
function. It is used to measure the impact of the dual network.  The
other two baselines are supervised learning (SL) frameworks. The first
supervised learning framework, \emph{Na\"ive}, uses a loss function
that minimizes the distance between the base case generation dispatch
estimates and the ground truth. The second supervised learning
framework is the \emph{Lagrangian Duality} (LD) framework. The
supervised learning models are impractical in this context, since they
require solving a very large number of instances optimally: they are
used to show that self-supervised learning is indeed the technique of
choice in this setting.

\begin{table}[!t]
\centering
\small
\setlength{\tabcolsep}{3pt}
\begin{tabular}{@{}lcccccc@{}}
\toprule
& \multicolumn{3}{c}{Before} & \multicolumn{3}{c}{After}\\
Test Case  & \#CV & \#BV & \#Cnst & \#CV & \#BV & \#Cnst \\
\midrule
\texttt{1354\_peg} & 3.3m   & 50.0k  & 6.7m   & 642.3k   & 50.0k    & 642.3k\\
\texttt{1888\_rte} & 4.8m   & 83.8k  & 9.7m   & 503.4k   & 83.8k    & 649.7k\\
\texttt{3022\_goc} & 14.6m  & 106.6k & 29.4m  & -         & -         & - \\
\texttt{4917\_goc} & 38.2m  & 321.5k & 77.1m  & -         & -         & - \\
\texttt{6515\_rte} & 64.9m  & 449.4k & 130.7m & -         & -         & - \\
\bottomrule
\end{tabular}
\caption{The Number of Binary and Continuous Variables Denoted as \#BV
  and \#CV and Constraints (\#Cnst) in Extensive SCOPF Problem Before
  and After Presolve. `$-$' indicates that Presolve runs of out of
  memory. $k$ and $m$ signify $10^3$ and $10^6$.}
\label{tab:numvar}
\end{table}

\begin{table}[!t]
\centering
\small
\setlength{\tabcolsep}{4pt}
\begin{tabular}{@{}lcccccc@{}}
\toprule
           & \multicolumn{3}{c}{Solving Time (s)} & \multicolumn{3}{c}{\#Iteration} \\
              \cmidrule(lr){2-4}\cmidrule(lr){5-7}
Test Case  & min & mean & max         & min & mean & max \\
\midrule
\texttt{1354\_peg}     & 42.66   & 97.10 & 1125.91 & 3 & 4.79 & 7 \\
\texttt{1888\_rte}     & 171.92  & 654.39 & 2565.6 & 3 & 4.93 & 6 \\
\texttt{3022\_goc}     & 628.67  & 6932.20 & 15844.52 & 8 & 10.71 & 17 \\
\texttt{4917\_goc}     & 3720.33 & 8035.86 & 15316.95 & 7 & 10.98 & 16 \\
\texttt{6515\_rte}     & 4648.31 & 9560.78 & 92430.47 & 4 & 6.65 & 10 \\
\bottomrule
\end{tabular}
\caption{Elapsed Time and Iterations for Solving SCOPF Test Instances Using CCGA (from \citep{PDLSCOPF}).}
\label{tab:CCGA_time}
\end{table}

\paragraph{Architectural Details}

Both the primal and dual networks consist of four fully-connected
layers, each followed by Rectified Linear Unit (ReLU) activations.
Layer normalization \citep{ba2016layer} is applied before the
fully connected layers for the primal network only.  The number of
hidden nodes in each fully connected layer is proportional to the
dimension of the input parameter, and is set to be $1.5\dim(\x)$.

The training uses a mini-batch size of 8 and a maximum of 1,000 epochs
for supervised baselines. PDL-SCOPF uses T=20 for the number of outer
iterations. Each iteration uses 2,000 forward/backward passes with a
minibatch size of 8 for both the primal and dual networks, resulting
in 80,000 iterations.  Note that, each iteration of PDF-SCOPF
generates instances {\em on the fly} from the distribution ${\cal P}$
for a total of $8\!\times\!80,000$ instances.  The implementation uses
PyTorch, and all models were trained on a computer with a NVIDIA Tesla
V100 GPU and an Intel Xeon 2.7GHz CPU. Averaged performance results
based on five independent training processes with different seeds are
reported.

\paragraph{Numerical Results}

\begin{table}[!t]
\centering
\small
\begin{tabular}{@{}l|cccc@{}}
\toprule
\begin{tabular}{@{}c@{}}\\ Test Case\end{tabular}
& \begin{tabular}{@{}c@{}}Na\"ive \\ (SL)\end{tabular} 
& \begin{tabular}{@{}c@{}}LD \\ (SL)\end{tabular} 
& \begin{tabular}{@{}c@{}} Penalty \\ (SSL)\end{tabular} 
& \begin{tabular}{@{}c@{}}PDL-SCOPF \\ (SSL) \end{tabular}\\
\midrule
\texttt{1354\_peg} & 0.022& 0.003 & 0.000 & 0.000 \\
\texttt{1888\_rte} & 0.044& 0.003 & 0.000 & 0.000 \\
\texttt{3022\_goc} & 0.088& 0.006 & 0.000 & 0.000 \\
\texttt{4917\_goc} & 0.067& 0.002 & 0.000 & 0.000 \\
\texttt{6515\_rte} & 0.001& 0.000 & 0.000 & 0.000 \\
\bottomrule       
\end{tabular}
\caption{Maximum Violations on Power Balance Constraints for Generator Contingency (in p.u.) (from \citep{PDLSCOPF}).}
\label{tab:pb_kg}
\end{table}

\begin{table}[!t]
\centering
\small
\setlength{\tabcolsep}{3pt}
\begin{tabular}{@{}l|r@{}lr@{}lr@{}lr@{}l@{}}
\toprule
\begin{tabular}{@{}c@{}}\\ Test Case\end{tabular}
& \multicolumn{2}{c}{ \begin{tabular}{@{}c@{}} Na\"ive \\ (SL)\end{tabular} } 
& \multicolumn{2}{c}{ \begin{tabular}{@{}c@{}} LD \\ (SL)\end{tabular} }
& \multicolumn{2}{c}{ \begin{tabular}{@{}c@{}} Penalty \\ (SSL)\end{tabular} } 
& \multicolumn{2}{c}{ \begin{tabular}{@{}c@{}} PDL-SCOPF \\ (SSL) \end{tabular} }\\
\midrule
\texttt{1354\_peg} & 14&(7.25\%) & 4&(2.07\%)  & 0&(0.00\%) & 0&(0.00\%) \\
\texttt{1888\_rte} & 31&(10.68\%)& 7&(2.41\%)  & 0&(0.00\%) & 0&(0.00\%) \\
\texttt{3022\_goc} & 57&(17.43\%)& 13&(3.98\%) & 0&(0.00\%) & 0&(0.00\%) \\
\texttt{4917\_goc} & 42&(7.41\%) & 12&(2.12\%) & 0&(0.00\%) & 0&(0.00\%) \\
\texttt{6515\_rte} & 2&(0.30\%)  & 0&(0.00\%)  & 0&(0.00\%) & 0&(0.00\%) \\
\bottomrule       
\end{tabular}
\caption{The Number of Generator Contingencies with Violated Power Balance Constraints (Percentages in Parenthesis) (from \citep{PDLSCOPF}).}
\label{tab:pb_kg_freq}
\end{table}

The numerical results compare PDL-SCOPF with the ground truth and the
baselines by evaluating the optimality gap and constraint violations.
Tables~\ref{tab:pb_kg}--\ref{tab:pb_kg_freq} show the violations of
the power balance constraint for generator
contingencies~\eqref{eq:scopf_ext_cnst_pb_kg}.  PDL-SCOPF and the SSL
penalty method have negligible violations of the power balance
constraints on all instances (they are within a tolerance of
$1\mathrm{e}{\textrm{-}4}$), thus producing a primal optimization
proxy. Table~\ref{tab:optgap} reports the mean optimality gap in
percentage, providing a comparison between the CCGA algorithm and its
learning counterparts. The results show that PDL-SCOPF is always the
strongest method with optimality gaps often below 1\% on these test
cases. PDL-SCOPF significantly dominates the Penalty (SSL) method,
showing the benefits of the dual network.

\begin{table}[!t]
\centering
\small
\begin{tabular}{l|cccc}
\toprule
\begin{tabular}{@{}c@{}}\\ Test Case\end{tabular}
& \begin{tabular}{@{}c@{}}Na\"ive\\ (SL)\end{tabular} 
& \begin{tabular}{@{}c@{}}LD\\ (SL)\end{tabular} 
& \begin{tabular}{@{}c@{}} Penalty \\ (SSL)\end{tabular} 
& \begin{tabular}{@{}c@{}}PDL-SCOPF \\ (SSL) \end{tabular}\\
\midrule
\texttt{1354\_peg}     & 13.447  & 2.700 & 2.533 & \textbf{0.856} \\
\texttt{1888\_rte}     & 22.115  & 2.436 & 4.969 & \textbf{1.960} \\
\texttt{3022\_goc}     & 159.116 & 8.008 & 1.312 & \textbf{0.983} \\ 
\texttt{4917\_goc}     & 47.212  & 2.096 & 0.454 & \textbf{0.210} \\
\texttt{6515\_rte}     & 2.419   & 1.292 & 2.069 & \textbf{0.815} \\
\bottomrule       
\end{tabular}
\caption{Mean Optimality Gap (\%) (best values in bold) (from \citep{PDLSCOPF}).}
\label{tab:optgap}
\end{table}

\begin{table}[!t]
\centering
\small
\begin{tabular}{l|rr}
\toprule
Test Case & Training & Sampling\\
\midrule
\texttt{1354\_peg} & 36min & 269hr 43min\\
\texttt{1888\_rte} & 43min & 1817hr 45min\\
\texttt{3022\_goc} & 1hr \phantom{0}7min & 19256hr \phantom{0}7min\\
\texttt{4917\_goc} & 1hr 59min & 22321hr 50min\\
\texttt{6515\_rte} & 2hr 54min & 26557hr 43min\\
\bottomrule       
\end{tabular}
\caption{Averaged Training Time (in GPU) for PDL-SCOPF and Accumulated CPU Time to Prepare the Supervised Training Dataset (from \citep{PDLSCOPF}).}
\label{tab:time}
\end{table}

\begin{table}[!t]
\centering
\small
\setlength{\tabcolsep}{3pt}
\begin{tabular}{l|cccc|c}
\toprule
            & \multicolumn{4}{c|}{PDL-SCOPF Inference Time (ms)} & \multirow{2}{*}{Speedup}\\
Test Case  & 8@GPU   &  1@GPU   &  8@CPU     & 1@CPU   & \\
\midrule                                                         
\texttt{1354\_peg} & 6.917   &  5.146   &  290.376   & 25.745  & 3771.60$\times$ \\ 
\texttt{1888\_rte} & 7.956   &  5.145   &  455.846   & 46.876  & 13960.02$\times$ \\
\texttt{3022\_goc} & 13.486  &  5.823   &  1484.019  & 183.879 & 37699.81$\times$\\
\texttt{4917\_goc} & 31.775  &  8.043   &  3975.468  & 479.881 & 16745.53$\times$\\
\texttt{6515\_rte} & 51.9370 &  10.576  &  6767.474  & 823.016 & 11616.76$\times$\\
\bottomrule       
\end{tabular}
\caption{
Averaged Inference Time of PDL-SCOPF with 1 Instance or 8 Instances on CPU or GPU in Milliseconds and Averaged Speedup over Gurobi on 1 CPU (from \citep{PDLSCOPF}).}
\label{tab:inference_time}
\end{table}

Table~\ref{tab:time} reports the elapsed GPU times for the PDL-SCOPF
training and the accumulated CPU time for generating instances for the
supervised baselines.  The time to generate solutions offline is
extremely high, even with CCGA.  The self-supervised learning
including PDL-SCOPF does not need the ground truth, which brings
substantial benefits in training time. Most interestingly, the
training time of 2 hours and 54 minutes for the French transmission
system (i.e., the largest test case \texttt{6515\_rte}) is
significantly smaller than the time required to solve the worst-case
SCOPF instance. Table~\ref{tab:inference_time} reports the inference
time in milliseconds of PDL-SCOPF on CPUs and GPUs for a single
instance or a batch of 8 instances. Even for the largest test case,
PDL-SCOPF provides a high-quality approximation to a single instance
in about 10 milliseconds and to a batch of instances in about 50
milliseconds. These results show that PDL-SCOPF is 4 orders of
magnitude faster than Gurobi on these instances.

{\em PDL-SCOPF highlights the step change that optimization learning may
bring.  Primal-Dual Learning can produce near-optimal solutions to a
complex power system applications in milliseconds, i.e., four orders
of magnitude faster than the best commercial solver and the state of
the art column and constraint generation algorithm. Opimization learning
thus expands the realm of applications which can benefit from optimization
technology.}

%% file: conclusion.tex
\section{Conclusion}
\label{sec:conclusion}

This paper presented the concept of {\em optimization learning}, a
methodology to design differentiable programs that can learn the
input/output mapping of parametric optimization problems. These
optimization proxies are trustworthy by design: they compute feasible
solutions to the underlying optimization problems, provide quality
guarantees on the returned solutions, and scale to large
instances. Optimization proxies combine traditional deep learning
technology with repair or completion layers to produce feasible
solutions. The paper also showed that optimization proxies can be
trained end-to-end in a self-supervised way.

The paper illustrated the impact of optimization learning on two
applications of interest to the power industry: the real-time risk
assessment of a transmission systems and the security-constrained
optimal power under N-1 Generator and Line Contingencies. In each
case, optimization proxies brings orders of magnitude improvements in
efficiency, which makes it possible to solve the applications in real
time with high accuracy, an outcome which that could not have been
achieved by state-of-the-art optimization technology.

There are many open issues in optimization learning. They include (1)
understanding how to derive effective repair layers for a wide range
of applications; (2) studying how to apply optimization learning
for combinatorial optimization problems, where the gradients are
typically not meaningful, and (3) applying Primal-Dual learning
to a variety of applications in nonlinear optimization.

Optimization learning is only one direction to fuse machine learning
and optimization. Learning to optimize, where machine learning is
used to speed up an existing algorithm is another avenue to leverage
the strenghts of both approaches.